\documentclass{birkjour}

 \newtheorem{thm}{Theorem}[section]
 \newtheorem{cor}[thm]{Corollary}
 \newtheorem{lem}[thm]{Lemma}
 \newtheorem{prop}[thm]{Proposition}
 \theoremstyle{definition}
 \newtheorem{defn}[thm]{Definition}
 \theoremstyle{remark}
 \newtheorem{rem}[thm]{Remark}
 
 \numberwithin{equation}{section}
\usepackage{amssymb}
\usepackage{amsmath}
\usepackage{graphicx}
\usepackage{amsmath}
\usepackage{txfonts}
\usepackage{layout}
\usepackage{latexsym}
\usepackage{float}
\usepackage{caption}
\usepackage{fullpage}
\usepackage{setspace}
\usepackage{color}
\usepackage[english]{babel}
\usepackage{dcpic, pictex}
\def\capa{\mathop{\hbox{\rm cap}}}

\def\R{\mathbb R}

\def\N{\mathbb N}
\def\N{\Bbb N}
\newcommand{\dx}{\:{\rm d}x}

\begin{document}

%
%
%
%
%
%
%
%
%

\title[ $(p,q)$-Laplacian double  obstacle problems]
 {Convergence results for the solutions of $(p,q)$-Laplacian double  obstacle problems on  irregular domains }

\author[R. Capitanelli ]{Raffaela Capitanelli}

\address{Dipartimento di Scienze di Base e Applicate per l'Ingegneria, $\lq\lq$ Sapienza" Università di Roma\\ 
Via A. Scarpa 16, 00161, Roma, Italy}

\email{raffaela.capitanelli@uniroma1.it}

\thanks{This work was completed with the support of our
\TeX-pert.}
\author{Salvatore Fragapane}
\address{$\lq\lq$Sapienza" Università di Roma}
\email{salvatore.fragapane@uniroma1.it}



\begin{abstract}
In this paper we study  double obstacle problems involving $(p,q)-$Laplace type operators. 
 In particular, we analyze the asymptotics of the solutions  on fractal and pre-fractal boundary domains.
\end{abstract}

\maketitle

\section{Introduction} \label{Intro} \setcounter{section}{1} \setcounter{equation}{0}

In this work we deal with double obstacle problems involving a family of $(p,q)-$Laplace type operators in fractal and pre-fractal boundary domains in $\R^2$.

Our motivation is due to the fact that  $(p,q)-$Laplace type operators  and irregular domains are useful  for  the study of many real problems.
In fact, $p$-Laplace operator  appears in the study of many concrete problems as non-Newtonian fluid mechanics, flows through porous media  (see \cite{DIA} and the references therein). The $(p, q)$-Laplacian has a wide range of applications in physical and related sciences, e.g. in biophysics,  quantum physics, plasma physics, solid state physics, chemical reaction design, and reaction-diffusion systems (see, for example,  \cite{BCM},  \cite{BPP},   \cite{BJ},  \cite{CI}, \cite{FOP}, \cite{MM}, \cite{Mar}, \cite{Z2} and the reference therein).
The limit operator (the $\infty$-Laplacian) also plays a leading role in problems like the mass transport (see \cite{MRT} and \cite{V} ) and the torsion creep (see, f.i., \cite{BDM} and \cite{K2}). 
Moreover, obstacle problems  appear in many other contexts: fluid filtration in porous media,
elasto-plasticity, optimal control, and financial math (see, for example, \cite{CAF} and \cite{FR}.)
With regard to fractals, what we mainly want to underline is their capacity to describe natural objects in a better way; in particular, they can represent the irregular  structure of several media (coasts, elements of the human body, etc.) more appropriately with respect to the \lq\lq classical  smooth" structures and therefore they find application in the modelization of various phenomena (see, for example, \cite{MAN} and \cite{MAN2}).\\

In this paper, for $ p>q$, $q\in[2,\infty)$ and $k\in\R$, we consider the following double obstacle problem on a fractal domain $\Omega_\alpha$ (see, for the definition, Section \ref{FpFB})
\begin{equation*}\label{PMpq} 
\min_{v\in\mathcal{H}_p}J_{p,q}(v),  \,\tag{ $\mathcal{PM}_{p,q}$}
\end{equation*}
with
\begin{equation}\label{Opi}
J_{p,q}(v)=\frac{1}{p}\int_{\Omega_\alpha}(k^2+|\nabla v|^2)^\frac{p}{2}\dx+\frac{1}{q}\int_{\Omega_\alpha}(k^2+|\nabla v|^2)^\frac{q}{2}\dx-\int_{\Omega_\alpha}fv\dx,
\end{equation}
\begin{equation}\label{convy}
\mathcal{H}_p=\{v\in W_g^{1,p}(\Omega_\alpha):\, \varphi_1\leq v\leq\varphi_2\text{ in }\Omega_\alpha\},
\end{equation}
and where $f\in L^{1}(\Omega_\alpha)$, $g\in W^{1,\infty}(\Omega_\alpha)$, $\varphi_1,\varphi_2\in C(\overline{\Omega}_{\alpha})$ are given
($W^{1, p}_g(\Omega):=\{v\in W^{1,p}(\Omega)\,:\, u=g \text{ on }\partial\Omega\}$).

As in \cite{BJ}, where the authors analyzed the homogeneous case without obstacles, we study the asymptotic behavior of the solutions, $u_{p,q}$, as $p\to\infty$, showing that, along subsequences, they converge to  solution of the following problem
\begin{equation*}\label{PMq} 
\min_{v\in\mathcal{H}}J_{q}(v), \text{ with } J_{q}(v)=\frac{1}{q}\int_{\Omega_\alpha}(k^2+|\nabla v|^2)^\frac{q}{2}\dx-\int_{\Omega_\alpha}fv\dx, \text{ if } L^2+k^2\leq1,\,\,\tag{ $\mathcal{PM}_{q}$}
\end{equation*}
\begin{equation*}\label{PMqL} 
\min_{v\in\mathcal{H}}||\nabla v||_\infty, \text{ if } L^2+k^2>1\,\,\tag{ $\mathcal{PM}_{q,L}$}
\end{equation*}
where $L$ is the Lipschitz constant of $g$ and
\begin{equation}\label{convyF}
\mathcal{H}=\Big\{v\in W_g^{1,\infty}(\Omega_\alpha):\, \varphi_1\leq v\leq\varphi_2\text{ in }\Omega_\alpha,\,||\nabla v||_\infty\leq\max\{1, \sqrt{L^2+k^2}\}\Big\},
\end{equation} (see Theorem \ref{AsymPtoI}).
We recall that asymptotic results for $p$-Laplacian have been studied   when  $p$ goes to $\infty,$ for example,  in \cite{BDM}, \cite{CV2}, \cite{IL} and \cite{MRT}.

Since the previous problems are defined on the domain $\Omega_\alpha$, which can be seen as the limit of appropriate domains $\Omega_\alpha^n$ (see Section \ref{FpFB}), it becomes natural to consider the corresponding approximating problems, that is the problems on the approximating domains $\Omega_\alpha^n$. For these problems, it is possible to prove analogous convergence results, as $p\to\infty$. Anyway, the introduction of the approximating problems opens another issue concerning the behavior of the solutions with respect to $n$, namely if it is possible to obtain a solution to Problems (\ref{PMpq}) and (\ref{PMq}) or (\ref{PMqL}) as limit, 
with respect to $n$, of the solutions of such problems.\\
The behavior of the solutions with respect to $n$ has been studied by many authors (see, f.i., \cite{AST}, \cite{C}, \cite{CV3}, \cite{CV1}, \cite{LV} and \cite{MV}). Nevertheless, as far as we know, in \cite{CF}, considering a double obstacle problem, the authors analyzed simultaneously both the behavior with respect to $p$ and $n$ for the first time. In fact, their analysis raises a question about the possibility of changing the order of the limit obtaining the same limit solution. Unfortunately, the lack of a uniqueness result for the case $p=\infty$ does not allow to the authors to affirm it (see \cite{ACJ}, \cite{BDM}, \cite{J2} and the references quoted in them). A first step in this direction is done in \cite{CV2}, where the authors give sufficient conditions which allow to obtain the convergence of the whole sequence. Then, in \cite{F} these conditions are used to state uniqueness results (in the case $p=\infty)$ for the same unilateral obstacle problem. 

The purpose of this paper is to give a complete answer to the asymptotic behavior of the solutions to the problems considered. Indeed, after the analysis of the behavior with respect to $p$, we will study the one with respect to $n$, showing that analogous results to the ones stated in \cite{CF} and \cite{F} hold even in the case here examined.\\
We stress the fact that in the proof of the convergence as $n\to \infty$ a necessary step is the construction of a sequence of functions, belonging to the approximating convex ($\mathcal{H}_{p,n}$ or $\mathcal{H}_n$), which converges to an element chosen in the corresponding final convex ($\mathcal{H}_{p}$ or $\mathcal{H}$). The hard part in this construction is to obtain functions satisfying all the conditions required from the convex, and this is due also to the irregular nature of the domains considered.
In this context, we emphasize how the introduction of a suitable coefficient functions and an integrability result for the gradient of the solutions will be fundamental tools in order to state our results. More precisely, for $p$ fixed and finite, the coefficient functions  allow to obtain a sequence of functions which converge to a solution of the problem on $\Omega_\alpha$, preserving some property, as long as its gradient has a greater summability with respect to the one of the natural space in which we search solutions. In this framework, the summability results we quoted play a crucial role; in particular, the approach and the techniques are the ones of \cite{GM}, \cite{GM2} and \cite{KK}.\\ 
For reasons of completeness, we specify that, beyond a greater summability of the gradient, further regularity results for the solution are present in the literature. In particular, regularity results for $p$-Laplacian obstacle problems in pre-fractal domains has been given in  \cite{CFV} (see also the references quoted there) and, in the same paper, these results are applied to give a optimal error estimate for the corresponding FEM problem following the approach used in \cite{G}.\\
Finally, we point out that it is possible to extend the  results of the present paper to other domains possibly with prefractal  and fractal boundaries 
if  these domains are \lq\lq Sobolev admissible domains" (see \cite{Dek}, \cite{Hinz}).

The plan of the paper is the following. In Section 2 we  introduce  the construction of fractal and pre-fractal boundary domains. In Section 3 the problem is introduced and  a corresponding integrability  result  for the gradient of the solution is obtained in Section 4. Section 5  is devoted to the asymptotic analysis with respect to $p$  and  Section 6 is devoted to the asymptotic analysis with respect to $n$. 

\section{Fractal and pre-fractal boundary domains} \label{FpFB} \setcounter{section}{2} 

In order to introduce the domains $\Omega^n_{\alpha}$ and $\Omega_{\alpha}$, which are the ones we will use in this paper, we need to remind how the construction of the Koch curve and the corresponding approximate pre-fractal curves works. So, let us recall the procedure which allows to obtain the n-th pre-fractal $K^n_{\alpha}$, $n\in\mathbb{N}$, of the Koch curves (see \cite{HU} for details and proofs).
Let us start considering, for instance, the line segment $K^0$ with endpoints $A(0,0)$ and $B(1,0)$ and let us introduce a family of four contractive similitudes $\Psi_{\alpha}=\{\psi_{1,\alpha},\dots,\psi_{4,\alpha}\}$  having $\alpha^{-1}$ as contraction factor, with $2<\alpha<4$, defined as follow:
$$
\psi_{1,\alpha}(z)=\frac{z}{\alpha},\,\quad\quad
\psi_{2,\alpha}(z)=\frac{z}{\alpha}e^{i\theta(\alpha)}+\frac{1}{\alpha},\qquad\quad$$
$$\psi_{3,\alpha}(z)=\frac{z}{\alpha}e^{-i\theta(\alpha)}+\frac12+i\sqrt{\frac{1}{\alpha}-\frac{1}{4}},\quad\quad
\psi_{4,\alpha}(z)=\frac{z-1}{\alpha}+1,
$$
with
\begin{equation}
\label{teta}
\theta(\alpha)=\arcsin\left( \frac{\sqrt{\alpha(4-\alpha)}}{2}\right).
\end{equation}
The first iteration makes us get a polygonal of four line segments. The following Figure \ref{Step1} shows this first step.

\begin{figure}[H]
	\begin{center}
		\includegraphics[height=2.5cm]{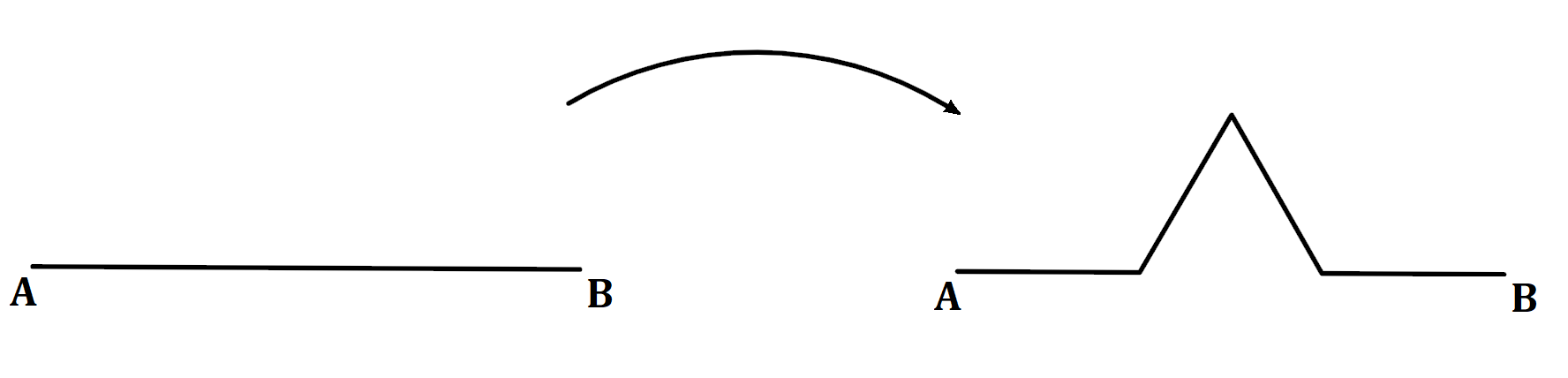}
	\end{center}
\caption{The starting segment $AB$ and the polygonal which come out of the first iteration in the case of $\alpha=3$.}
\label{Step1}
\end{figure}

\begin{figure}[H]
	\begin{center}
		\includegraphics[height=3.5cm]{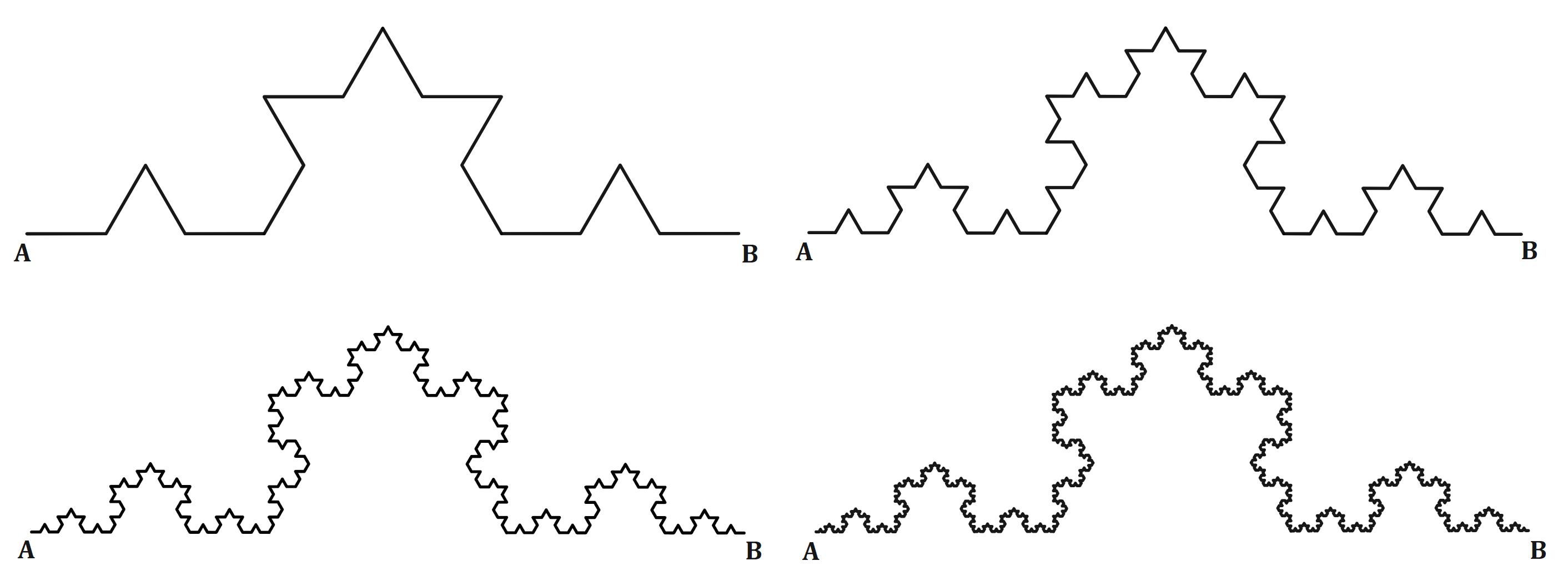}
	\end{center}
\caption{$K^n_3$ for $n=2$, $n=3$, $n=4$ and $K_3$.}
\label{KochN}
\end{figure}

In general, at every step each segment of the polygonal is replaced with a rescaled copy of the one in the basic step (see Figure \ref{Step1}).\\
Then, for each $n\in\N$, we set
$$ K_\alpha^{n}=\bigcup_{i=1}^4 \psi_{i,\alpha}(K_\alpha^{n-1})=\bigcup_{i|n} K_{\alpha}^{i|n}, \text{ with } K_{\alpha}^{i|n}=\psi_{i|n,\alpha}(K^0),$$
where $\psi_{i|n,\alpha}=\psi_{i_1,\alpha}\circ\psi_{i_2,\alpha}\circ\cdots\circ\psi_{i_n,\alpha}$ is the map associated with an arbitrary $n-$tuple of indices $i|n=(i_1, i_2, \ldots, i_n)\in \{1,\dots,4\}^n$, for each integer $n>0$ and $\psi_{i|n,\alpha}=id$ in $\R^2$ if $n=0$.\\
Figure \ref{KochN} shows the result of some steps of the procedure just recalled and the final curve.

As $n\to\infty$, we get that the curves $K_{\alpha}^n$ converge to the fractal curve $K_\alpha$ in the Hausdorff metric. Moreover, $K_\alpha$ is the unique compact set which is invariant on $\Psi_{\alpha}$ and  $d_f=\frac{\ln4}{\ln\alpha}$ is its Hausdorff dimension.\\

\begin{figure}[H]
	\begin{center}
		\includegraphics[height=7.5cm]{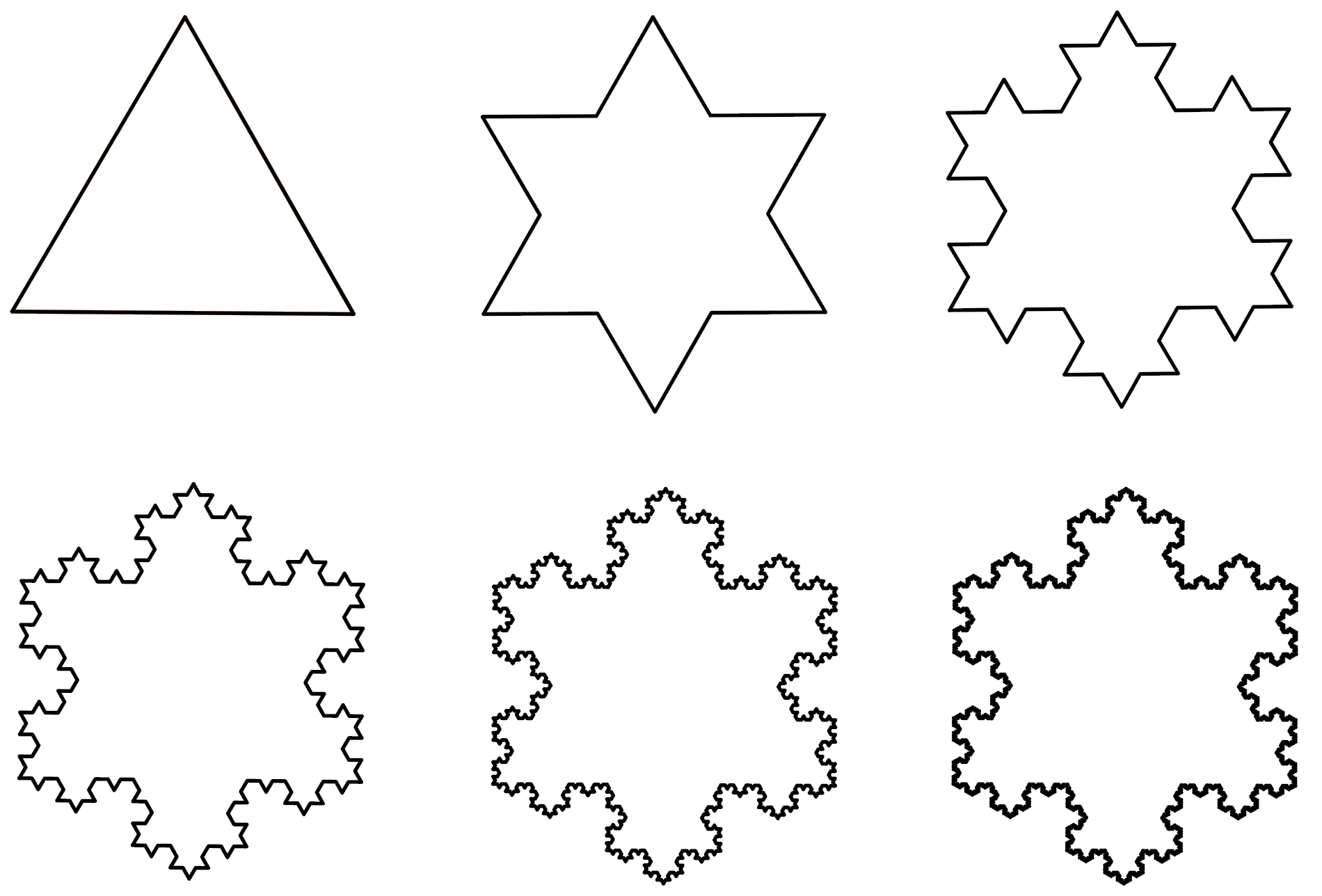}
	\end{center}
\caption{$\Omega^0$, $\Omega^n_3$ for $n=1$, $n=2$, $n=3$, $n=4$ and $\Omega_3$.}
\label{Omega}
\end{figure}

Now, we denote with $\Omega^n_\alpha$ the domain obtained starting to any regular polygon $\Omega^0$ (triangle, square, etc.) and replacing each of its sides with the n-th pre-fractal Koch curve $K_{\alpha}^n$. These domains are non-convex, polygonal, with an increasing number of sides and, at the limit, they develop the fractal geometry of $\Omega_{\alpha}$, i.e. the domain having a fractal boundary formed by the union of Koch curves 
(see Figure \ref{Omega}).

\section{Double obstacle problem}
\setcounter{section}{3} 

Let $p>q$ be, with $q\in[2,\infty)$ fixed. Given $f\in L^{1}(\Omega_\alpha)$, $g\in W^{1,\infty}(\Omega_\alpha)$ and $\varphi_1,\varphi_2\in C(\overline{\Omega}_{\alpha})$, 
Problem (\ref{PMpq}) is equivalent to the following variational inequality 
\begin{equation*}\label{Ppq} 
\text{find } u_{p,q}\in\mathcal{H}_p\,\,:\,\,a_{p}(u_{p,q},v-u_{p,q})+a_{q}(u_{p,q},v-u_{p,q})-\int_{\Omega_\alpha}f(v-u_{p,q})\dx \geqslant 0,\forall v\in \mathcal{H}_p,\,\tag{ $\mathcal{P}_{p,q}$}
\end{equation*}
where
\begin{equation}\label{ops}
a_p(u,v)=\int_{\Omega_\alpha}(k^2+|\nabla u|^2)^{\frac{p-2}{2}}\nabla u\nabla v\dx,
\end{equation}
and $\mathcal{H}_p$ is defined in (\ref{convy}).\\
Moreover (see, for instance, \cite{T}), if  $\mathcal{H}_{p}$ is non-empty,  as  the functional $J_{p,q}(v)$ is convex, weakly lower semi-continuous and coercive, then Problem (\ref{PMpq}) has a minimizer $u_{p,q}$ in $\mathcal{H}_{p}$.

\begin{prop}\label{UniqPpq}
Problem (\ref{Ppq})  admits unique solution.
\end{prop}
\proof Let $u_1$ and $u_2$ solutions to Problem (\ref{Ppq}). It holds that
$$a_{p}(u_1,u_2-u_1)+a_{q}(u_1,u_2-u_1)-\int_{\Omega_\alpha}f(u_2-u_1)\dx \geq 0$$
and
$$-a_{p}(u_2,u_2-u_1)-a_{q}(u_2,u_2-u_1)+\int_{\Omega_\alpha}f(u_2-u_1)\dx \geq 0.$$
So, by the definition of $a_p(u,v)$ and summing previous relations, we get
$$\int_{\Omega_\alpha}(k^2+|\nabla u_1|^2)^{\frac{p-2}{2}}\nabla u_1\nabla(u_2-u_1)\dx-\int_{\Omega_\alpha}(k^2+|\nabla u_2|^2)^{\frac{p-2}{2}}\nabla u_2\nabla(u_2-u_1)\dx+$$
$$+\int_{\Omega_\alpha}(k^2+|\nabla u_1|^2)^{\frac{q-2}{2}}\nabla u_1\nabla(u_2-u_1)\dx-\int_{\Omega_\alpha}(k^2+|\nabla u_2|^2)^{\frac{q-2}{2}}\nabla u_2\nabla(u_2-u_1)\dx\geq0,$$
that is
$$\int_{\Omega_\alpha}\Big[(k^2+|\nabla u_2|^2)^{\frac{p-2}{2}}\nabla u_2\nabla(u_2-u_1)-(k^2+|\nabla u_1|^2)^{\frac{p-2}{2}}\nabla u_1\nabla(u_2-u_1)\Big]\dx+$$
$$+\int_{\Omega_\alpha}\Big[(k^2+|\nabla u_2|^2)^{\frac{q-2}{2}}\nabla u_2\nabla(u_2-u_1)-(k^2+|\nabla u_1|^2)^{\frac{q-2}{2}}\nabla u_1\nabla(u_2-u_1)\Big]\dx\leq0.$$
Thus, by Lemma 2.1 in \cite{LB}, we deduce that there exist $C=C(k)>0$ such that
$$C\Big(||\nabla(u_2-u_1)||_p+||\nabla(u_2-u_1)||_q\Big)\leq0\Longleftrightarrow||\nabla(u_2-u_1)||_p+||\nabla(u_2-u_1)||_q\leq0$$
and then
$$u_1(x)=u_2(x)+c,\,\,c\in\R$$
By taking $y\in\partial\Omega_\alpha$, we have $u_1(y)=u_2(y)=g(y)$ and then
$$g(y)=g(y)+c\Longleftrightarrow c=0.$$
Hence, $u_1=u_2$.
\endproof

\begin{rem}\label{equi}
Let us consider $u_{p,q}\in\mathcal{H}_p$ solution to the problem
\begin{equation}\label{Ppq2} 
\text{find } u\in\mathcal{H}_p\,\,:\,\,a_{p}(u,v-u)+a_{q}(u,v-u)-\int_{\Omega_\alpha}f(v-u)\dx \geqslant 0,\forall v\in \mathcal{K}_p.\,\,\tag{ $\mathcal{P}_{p,q}$}
\end{equation} 
Denoting 
$$\bar{u}_E=\frac{1}{|E|}\int_{E}u\dx,$$
with $E\subseteq\Omega_\alpha$, we have that $(\mathcal{P}_{p,q})$ is equivalent to
$$a_{p}(u_{p,q}-\bar{u}_{p,q,E},v-\bar{u}_{p,q,E}+\bar{u}_{p,q,E}-u_{p,q})+a_{q}(u_{p,q}-\bar{u}_{p,q,E},v-\bar{u}_{p,q,E}+\bar{u}_{p,q,E}-u_{p,q})+$$ $$-\int_{\Omega_\alpha}f(v-\bar{u}_{p,q,E}+\bar{u}_{p,q,E}-u_{p,q})\dx \geqslant 0,\forall v\in \mathcal{K}_p.$$
Hence, defining
$$\hat{v}=v-\bar{u}_{p,q,E},\qquad \hat{\varphi}_i=\varphi_i-\bar{u}_{p,q,E},\, i=1,2\,,\qquad\text{and}\qquad\hat{u}_{p,q}=u_{p,q}-\bar{u}_{p,q,E},$$
we get that $\hat{u}_{p,q}\in\mathcal{H}_{p,-\bar{u}_{p,q,E}}=\{\hat{v}\in W_{g-\bar{u}_{p,E}}^{1,p}(\Omega_\alpha):\,\hat{\varphi}_1\leq\hat{v}\leq\hat{\varphi}_2 \text{ in }\Omega_\alpha\}$ is solution to problem
\begin{equation}\label{Ppqu}
\text{find } w\in\mathcal{H}_{p,-\bar{u}_{p,q,E}}\,\,:\,\,a_{p}(w,\hat{v}-w)+a_{q}(w,\hat{v}-w)-\int_{\Omega_\alpha}f(\hat{v}-w)\dx \geqslant 0,\forall \hat{v}\in \mathcal{K}_{p,-\bar{u}_{p,q,E}}.
\end{equation}
\end{rem}

\section{Integrability result}\label{SeI}
\setcounter{section}{4} 

A first step to do in order to deal with the asymptotic behavior, as $n\to\infty$, consists in an integrability result.

 The issue of the integrability of the gradient was widely studied by many authors; here, in particular, we refer to the approach and the results of \cite{GM}, \cite{GM2} and \cite{KK}. On the one hand, in \cite{GM} and \cite{GM2} the authors gave a global integrability result for obstacle problems, assuming that the boundary was $p$-Poincaré thick. On the other hand, in \cite{KK} the authors obtained an analogous result for a  problem without obstacle and involving an operator satisfying the same assumptions, but requiring a weaker condition on the boundary. \\

In the following subsection, we simply recall  the  previous quoted  results making only the necessary changes or assumptions to adapt them to our case.

\subsection{Preliminary tools}\label{Pt}

Let us recall some definitions and crucial lemmas which are preparatory to the integrability result.\\
We point out that from now on, we adopt the following notation:
$$\cdot\,\,\,B_{\rho}(x)=\{x\in\R^2\,:\, ||x||<\rho\},\,\,\rho>0;$$
$$\cdot\,\,\,B_{\rho} \text{ a generic ball having radius equal to } \rho.$$

\begin{lem}\label{Poinc}(\textbf{Poincar\'e's inequalities})\\
Let $\Omega$ be a bounded open subset of $\R^2$, with diameter $d$, and let $p\geq1$ be.\\

(i) If $u\in W_0^{1,p}(\Omega)$, then
	\begin{equation}\label{Poinc1}
		||u||_p\leq \frac{d}{p^{\frac{1}{p}}}||\nabla   u||_p\leq d||\nabla u||_p,
	\end{equation}
(see Theorem 12.17 in \cite{L}).\\
	
(ii) If $\Omega$ is convex and $u\in W^{1,p}(\Omega)$, then $\exists C=C(p)>0$ such that
\begin{equation}\label{Poinc2}
||u-\bar{u}_{\Omega}||_{p,\Omega}\leq Cd||\nabla u||_{p,\Omega},
\end{equation}
(see Theorem 12.30 in \cite{L}).
\end{lem}

\begin{lem}\label{Sob}(\textbf{Sobolev-Poincar\'e's inequalities})\\
Let us consider $p\in[1,2)$ and $u\in W^{1,p}(\R^2)$, then there exists $C=C(p)>0$ such that	
\begin{equation}\label{Sob2}
||u-\bar{u}_{B_r}||_{p^*,B_r}\leq Cr||\nabla u||_{p,B_r}, \text{ with } p^*=\frac{2p}{2-p},
\end{equation}
for every $B_r\subset\R^2$ (see Theorems 3.16 and 3.20 in \cite{K}).
\end{lem}

For the following definitions and lemmas (and more informations about them) we refer, for instance, to \cite{HKM}, \cite{KK}, \cite{K}, \cite{Z} and the references quoted there.
\begin{defn}
	The $p$-capacity of a compact set $K\subset\Omega$ in $\Omega$ is defined as
	$${\capa}_p(K;\Omega)=\inf_{\substack{u\in C_0^{\infty}(\Omega)\\u=1 \text{ in } K}}\Big\{\int_{\Omega}|\nabla u|^p\dx\Big\}$$
	and for an arbitrary set $A\subset\Omega$ is defined as 
	$${\capa}_p(A;\Omega)=\inf_{\substack{A\subset E\subset\Omega\\ E\text{ open}}}\sup_{\substack{K\subset E\\K \text{ compact}}}\Big\{{\capa}_p(K;\Omega)\Big\},$$
\end{defn}

\begin{rem} Since we are in $\R^2$,
$${\capa}_p(\overline{B}_r;B_{2r})=Cr^{2-p},$$
	where $C=C(p)>0$.
\end{rem}

\begin{defn}
	We say that a function $u$ is $p$-quasicontinuous in $\Omega$ if, for any $\varepsilon>0$, there exists an open set $A$ with ${\capa}_p(A;\Omega)<\varepsilon$ such that $u$ restricted to $\Omega\setminus A$ is continuous.
\end{defn}

\begin{rem}
	For any $v\in W^{1,p}(\Omega)$, $p>1$, there exists a $p$-quasicontinuous fuction $w\in W^{1,p}(\Omega)$ such that $v=w$ q.e. in $\Omega$. Moreover, this representative $w$ is unique, in the sense that every other is equal to it except at most in a set of $p$-capacity equal to zero.
\end{rem}

\begin{defn}
	We say that a set $\Omega\subset\R^m,\,m\in\N$ is uniformly $p$-thick,
	$1<p<\infty$, if there exist positive constants $C$ and $\overline{r}$ such that
	$${\capa}_p(\Omega\cap\overline{B}_r(x);B_{2r}(x))\geq C{\capa}_p(\overline{B}_r(x);B_{2r}(x))$$
	whenever $x\in\Omega$ and $0<r<\overline{r}$.
\end{defn}

\begin{rem}\label{pthi}
	In our case this condition is not restrictive, since it is automatically satisfied; indeed if $p>m$ each non-empty set is uniform $p$-thick. Moreover, the complement of any proper simply connected subdomain of $\R^2$ is uniformly $p$-thick for all $p>1$.
\end{rem}
From now on, we denote
$$\cdot\,\,\,\fint_{\Omega}u\dx:=\frac{1}{|\Omega|}\int_{\Omega}u\dx,$$
where $\Omega\subset\R^2$, $u:\Omega\to\R$.
\begin{lem} \label{Som}
	Let $B\subset\R^2$ be a fixed ball and let $g,h\in L^s(B)$ be non-negative functions. If for some $s>1$ 
	$$\fint_{B_r}g^s\dx\leq C\Big(\fint_{B_{2r}}g\dx\Big)^{s}+\fint_{B_{2r}}h^s\dx$$	
	for each ball $B_r$, $r>0$, with $B_{2r}\subset B$,
	then there exists $\varepsilon=\varepsilon(C,s)>0$ such that for all $t\in[s,s+\varepsilon]$ it holds
	$$\Big(\fint_{B_{2r}}g^t\dx\Big)^{\frac{1}{t}}\leq C\Big(\fint_{B_{2r}}g^s\dx\Big)^{\frac{1}{s}}+\Big(\fint_{B_{2r}}h^t\dx\Big)^{\frac{1}{t}}$$	
\end{lem}	
\noindent(see Lemma 2.1 in \cite{KK}).

\begin{lem} \label{Cap}
Let $u$ be a $t$-quasicontinuous function in $W^{1,t}(B_r)$, where $t>1$ and $B_r\subset\R^2$, $r>0$.\\ 
Let $N(u)=\{x\in B_r: u(x)=0\}$. Then
$$\Big(\fint_{B_r}u^{\alpha t}\dx\Big)^{\frac{1}{\alpha t}}\leq C\Big(\frac{1}{{\capa}_t(N(u);B_{2r})}\int_{B_{r}}|\nabla u|^t\dx\Big)^{\frac{1}{t}},$$
where $C=C(t)>0$ and 
$$\alpha=\begin{cases}
\frac{2}{2-t}, \text{ if } 1<t<2\\
2, \text{ if } t\geq2
\end{cases} .$$
\end{lem}	
\noindent(see Lemma 3.1 in \cite{KK}).\\
This following last technical results is crucial in the proof of the integrability. We point out that it is an adaptation to our case of the analogous  ones proved in \cite{GM}, \cite{GM2} and \cite{KK}, and its proof uses the same arguments.
\begin{lem} \label{Ineq}
Let $p>q$, $q\in[2,\infty)$, $f\in L^{p'}(\Omega_\alpha)$, with $\frac{1}{p}+\frac{1}{p'}=1$, $\varphi_1,\varphi_2\in W^{1,p}(\Omega_\alpha)$ and $u$ solution to Problem (\ref{Ppq}).
Let $r>0$ such that $B_{r}\subset\Omega_\alpha$ and $\psi\in C_0^{\infty}(\Omega_\alpha)$, with $0\leq\psi\leq1$. Then there exists $C=C(k,p,q,\Omega_\alpha)>0$ such that:
$$\int_{\Omega_\alpha}|\psi|^p|\nabla u|\dx\leq C\Big[\int_{\Omega_\alpha}|\psi|^p(|\nabla\varphi_1|^p+|\nabla\varphi_2|^p)\dx+\int_{\Omega_\alpha}|\psi|^{\frac{p}{p-1}}\Big(|\nabla\varphi_1|^{\frac{p}{p-1}}+|\nabla\varphi_2|^{\frac{p}{p-1}}\Big)\dx+$$ $$+\int_{\Omega_\alpha}|\psi|^{\frac{p}{p-q+1}}\Big(|\nabla\varphi_1|^{\frac{p}{p-q+1}}+|\nabla\varphi_2|^{\frac{p}{p-q+1}}\Big)\dx+
\int_{\Omega_\alpha}|\nabla\psi|^p(|u-\bar{u}_{B_{r}}|^p+|\varphi_1-\overline{\varphi_1}_{B_{r}}|^p+$$ $$+|\varphi_2-\overline{\varphi_2}_{B_{r}}|^p)\dx+\int_{\Omega_\alpha}|\nabla\psi|^{\frac{p}{p-1}}(|u-\bar{u}_{B_{r}}|^{\frac{p}{p-1}}+|\varphi_1-\overline{\varphi_1}_{B_{r}}|^{\frac{p}{p-1}}+|\varphi_2-\overline{\varphi_2}_{B_{r}}|^{\frac{p}{p-1}})\dx+$$
\begin{equation}\label{Ineq1}
+\int_{\Omega_\alpha}|\nabla\psi|^{\frac{p}{p-q+1}}(|u-\bar{u}_{B_{r}}|^{\frac{p}{p-q+1}}+|\varphi_1-\overline{\varphi_1}_{B_{r}}|^{\frac{p}{p-q+1}}+|\varphi_2-\overline{\varphi_2}_{B_{r}}|^{\frac{p}{p-q+1}})\dx+\int_{\Omega_\alpha}|f|^{\frac{p}{p-1}}|\psi|^{\frac{p}{p-1}}\dx\Big];
\end{equation}
$$\int_{\Omega_\alpha}|\psi|^p|\nabla u|\dx\leq C\Big[\int_{\Omega_\alpha}|\psi|^p(|\nabla\varphi_1|^p
+|\nabla g|^p+|\nabla\varphi_2|^p)\dx+\int_{\Omega_\alpha}|\psi|^{\frac{p}{p-1}}\Big(|\nabla\varphi_1|^{\frac{p}{p-1}}+|\nabla g|^{\frac{p}{p-1}}+|\nabla\varphi_2|^{\frac{p}{p-1}}\Big)\dx+$$
$$+\int_{\Omega_\alpha}|\psi|^{\frac{p}{p-q+1}}\Big(|\nabla\varphi_1|^{\frac{p}{p-q+1}}+|\nabla g|^{\frac{p}{p-q+1}}+|\nabla\varphi_2|^{\frac{p}{p-q+1}}\Big)\dx+\int_{\Omega_\alpha}|\nabla\psi|^p|u-w|^p\dx+\int_{\Omega_\alpha}|\nabla\psi|^{\frac{p}{p-1}}|u-w|^{\frac{p}{p-1}}\dx+$$
\begin{equation}\label{Ineq2}	
+\int_{\Omega_\alpha}|\nabla\psi|^{\frac{p}{p-q+1}}|u-w|^{\frac{p}{p-q+1}}\dx+\int_{\Omega_\alpha}|f|^{\frac{p}{p-1}}|\psi|^{\frac{p}{p-1}}\dx\Big],
\end{equation}
where $w=(\varphi_1\vee g)\land\varphi_2$, with $g\in W^{1,\infty}(\Omega_\alpha)$.
\end{lem}	
\proof
Let us prove (\ref{Ineq1}). The proof is not difficult, but a little tangled. \\
Since $u\in\mathcal{H}_p$ is solution to Problem (\ref{Ppq}), thanks to Remark \ref{equi}, we have
\begin{equation}\label{EquiPp}
\int_{\Omega_\alpha}(k^2+|\nabla\hat{u}|^2)^{\frac{p-2}{2}}\nabla\hat{u}\nabla(\hat{v}-\hat{u})\dx
+\int_{\Omega_\alpha}(k^2+|\nabla\hat{u}|^2)^{\frac{q-2}{2}}\nabla\hat{u}\nabla(\hat{v}-\hat{u})\dx-\int_{\Omega_\alpha}f(\hat{v}-\hat{u})\dx \geqslant 0
\end{equation}
$\forall \hat{v}\in \mathcal{H}_{p,-\bar{u}_{B_r}},$
with $\hat{u}=u-\bar{u}_{B_{r}}$, $\hat{u}\in\mathcal{H}_{p,-\bar{u}_{B_{r}}}$.\\ 
Now, let us consider $$\hat{v}=u-\bar{u}_{B_{r}}-\psi^p(u-\bar{u}_{B_{r}})+\psi^pw=(1-\psi^p)(u-\bar{u}_{B_{r}})+\psi^pw,$$ 
with 
\begin{equation}\label{w}
w=(\varphi_1-\bar{u}_{B_r})^+\land(\varphi_2-\bar{u}_{B_r})^+-(\varphi_2-\bar{u}_{B_r})^-=
\begin{cases}
\varphi_1-\bar{u}_{B_r}, \text{ if } \varphi_1> \bar{u}_{B_r}\\
0, \text{ if } \varphi_1\leq \bar{u}_{B_r}\leq\varphi_2\\
\varphi_2-\bar{u}_{B_r}, \text{ if } \varphi_2< \bar{u}_{B_r}.
\end{cases}
\end{equation}
So, we deduce that $\hat{v}\in\mathcal{H}_{p,-\bar{u}_{B_r}}$.\\
Moreover, it holds
$$|w|\leq
\begin{cases}
|\varphi_1-\overline{\varphi_1}_{B_r}|, \text{ if } \varphi_2\geq\bar{u}_{B_r}\\
|\varphi_2-\overline{\varphi_2}_{B_r}|, \text{ if } \varphi_2< \bar{u}_{B_r}
\end{cases}$$
and then $|w|\leq |\varphi_1-\overline{\varphi_1}_{B_r}|+|\varphi_2-\overline{\varphi_2}_{B_r}|$.
Furthermore, by (\ref{w}), we deduce that $|\nabla w|\leq |\nabla\varphi_1|+|\nabla\varphi_2|$.\\
With this choice of $\hat{v}$ we have 
$$\nabla(\hat{v}-\hat{u})=\nabla(-\psi^p(u-\bar{u}_{B_{r}})+\psi^pw)=p\psi^{p-1}\nabla\psi (-(u-\bar{u}_{B_{r}})+w)-\psi^p\nabla u+\psi^p\nabla w.$$ 
Then, by (\ref{EquiPp}) and previous observations, we have:
$$0\leq\int_{\Omega_\alpha}(k^2+|\nabla u|^2)^{\frac{p-2}{2}}\nabla u\cdot[p\psi^{p-1}\nabla\psi (-(u-\bar{u}_{B_{r}})+w)-\psi^p\nabla u+\psi^p\nabla w]\dx+$$
$$+\int_{\Omega_\alpha}(k^2+|\nabla u|^2)^{\frac{q-2}{2}}\nabla u\cdot[p\psi^{p-1}\nabla\psi (-(u-\bar{u}_{B_{r}})+w)-\psi^p\nabla u+\psi^p\nabla w]\dx+$$
$$-
\int_{\Omega_\alpha}f(-\psi^p(u-\bar{u}_{B_{r}})+\psi^pw)\dx\leq$$

$$\leq-\int_{\Omega_\alpha}(k^2+|\nabla u|^2)^{\frac{p-2}{2}}|\nabla u|^2|\psi|^p\dx+\int_{\Omega_\alpha}(k^2+|\nabla u|^2)^{\frac{p-2}{2}}|\nabla u|(|\nabla\varphi_1|+|\nabla\varphi_2|)|\psi|^p\dx+$$
$$+p\int_{\Omega_\alpha}(k^2+|\nabla u|^2)^{\frac{p-2}{2}}|\nabla u||\nabla \psi||\psi|^{p-1}(|u-\bar{u}_{B_{r}}|+|\varphi_1-\overline{\varphi_1}_{B_{r}}|+|\varphi_2-\overline{\varphi_2}_{B_{r}}|)\dx+$$ $$-\int_{\Omega_\alpha}(k^2+|\nabla u|^2)^{\frac{q-2}{2}}|\nabla u|^2|\psi|^p\dx+$$
$$+\int_{\Omega_\alpha}(k^2+|\nabla u|^2)^{\frac{q-2}{2}}|\nabla u|(|\nabla\varphi_1|+|\nabla\varphi_2|)|\psi|^p\dx+$$ $$+p\int_{\Omega_\alpha}(k^2+|\nabla u|^2)^{\frac{q-2}{2}}|\nabla u||\nabla \psi||\psi|^{p-1}(|u-\bar{u}_{B_{r}}|+|\varphi_1-\overline{\varphi_1}_{B_{r}}|+|\varphi_2-\overline{\varphi_2}_{B_{r}}|)\dx+$$
$$+\int_{\Omega_\alpha}|f\psi||\psi|^{p-1}(|u-\bar{u}_{B_{r}}|+|\varphi_1-\overline{\varphi_1}_{B_{r}}|+|\varphi_2-\overline{\varphi_2}_{B_{r}}|)\dx.$$
Now, since
\begin{equation}\label{Pro}
\forall a,b\geq0,\,\,\forall r\geq1\,\,\exists\,\,C=C(r)>0\,\,:\,\,(a+b)^r\leq C(a^r+b^r),
\end{equation}
and by the fact that $|\psi|^a\leq|\psi|^b$, for any $a\geq b>0$, we get
$$0\leq-\int_{\Omega_\alpha}|\nabla u|^p|\psi|^p\dx+C_{k,p}\int_{\Omega_\alpha}|\nabla u|(|\nabla\varphi_1|+|\nabla\varphi_2|)|\psi|^2\dx+C_p\int_{\Omega_\alpha}|\nabla u|^{p-1}(|\nabla\varphi_1|+|\nabla\varphi_2|)|\psi|^p\dx+$$
$$+pC_{k,p}\int_{\Omega_\alpha}|\nabla u||\nabla \psi||\psi|(|u-\bar{u}_{B_{r}}|+|\varphi_1-\overline{\varphi_1}_{B_{r}}|+|\varphi_2-\overline{\varphi_2}_{B_{r}}|)\dx+pC_p\int_{\Omega_\alpha}|\nabla u|^{p-1}|\nabla \psi||\psi|^{p-1}(|u-\bar{u}_{B_{r}}|+|\varphi_1-\overline{\varphi_1}_{B_{r}}|+|\varphi_2-\overline{\varphi_2}_{B_{r}}|)\dx+$$
$$-\int_{\Omega_\alpha}|\nabla u|^q|\psi|^p\dx+C_{k,q}\int_{\Omega_\alpha}|\nabla u|(|\nabla\varphi_1|+|\nabla\varphi_2|)|\psi|^2\dx+C_q\int_{\Omega_\alpha}|\nabla u|^{q-1}(|\nabla\varphi_1|+|\nabla\varphi_2|)|\psi|^q\dx+$$
$$+pC_{k,q}\int_{\Omega_\alpha}|\nabla u||\nabla \psi||\psi|(|u-\bar{u}_{B_{r}}|+|\varphi_1-\overline{\varphi_1}_{B_{r}}|+|\varphi_2-\overline{\varphi_2}_{B_{r}}|)\dx+pC_q\int_{\Omega_\alpha}|\nabla u|^{q-1}|\nabla \psi||\psi|^{q-1}(|u-\bar{u}_{B_{r}}|+$$
$$+|\varphi_1-\overline{\varphi_1}_{B_{r}}|+|\varphi_2-\overline{\varphi_2}_{B_{r}}|)\dx+$$
$$+\int_{\Omega_\alpha}|f\psi||\psi|(|u-\bar{u}_{B_{r}}|+|\varphi_1-\overline{\varphi_1}_{B_{r}}|+|\varphi_2-\overline{\varphi_2}_{B_{r}}|)\dx.$$
Applying Young's inequality with conjugate exponents $p$ and $\frac{p}{p-1}$ or $\frac{p}{q-1}$ and $\frac{p}{p-q+1}$, we obtain:
$$0\leq-\int_{\Omega_\alpha}|\nabla u|^p|\psi|^p\dx+\frac{\sigma}{p}\int_{\Omega_\alpha}|\nabla u|^p|\psi|^p\dx+\frac{\overline{C}_{k,p}(p-1)\sigma^{-\frac{1}{p-1}}}{p}\int_{\Omega_\alpha}(|\nabla\varphi_1|^{\frac{p}{p-1}}+|\nabla\varphi_2|^{\frac{p}{p-1}})|\psi|^{\frac{p}{p-1}}\dx+$$
$$+\frac{\sigma_1(p-1)}{p}\int_{\Omega_\alpha}|\nabla u|^p|\psi|^p\dx+\frac{\sigma_1^{-(p-1)}}{p}\overline{C}_p\int_{\Omega_\alpha}(|\nabla\varphi_1|^p+|\nabla\varphi_2|^p)|\psi|^p\dx+\sigma_2\int_{\Omega_\alpha}|\nabla u|^p|\psi|^p\dx+$$
$$+\overline{C}_{k,p}(p-1)\sigma_2^{-\frac{1}{p-1}}\int_{\Omega_\alpha}|\nabla \psi|^{\frac{p}{p-1}}(|u-\bar{u}_{B_{r}}|^{\frac{p}{p-1}}+|\varphi_1-\overline{\varphi_1}_{B_{r}}|^{\frac{p}{p-1}}+|\varphi_2-\overline{\varphi_2}_{B_{r}}|^{\frac{p}{p-1}})\dx+$$ $$+\sigma_3(p-1)\int_{\Omega_\alpha}|\nabla u|^p|\psi|^p\dx+$$
$$+\overline{C}_{p}\sigma_3^{-(p-1)}\int_{\Omega_\alpha}|\nabla \psi|^p(|u-\bar{u}_{B_{r}}|^p+|\varphi_1-\overline{\varphi_1}_{B_{r}}|^p+|\varphi_2-\overline{\varphi_2}_{B_{r}}|^p)\dx-\int_{\Omega_\alpha}|\nabla u|^q|\psi|^q\dx+$$ $$ +\frac{\sigma_4}{p}\int_{\Omega_\alpha}|\nabla u|^p|\psi|^p\dx+$$
$$+\frac{C_{k,p,q}(p-1)\sigma_4^{-\frac{1}{p-1}}}{p}\int_{\Omega_\alpha}(|\nabla\varphi_1|^{\frac{p}{p-1}}+|\nabla\varphi_2|^{\frac{p}{p-1}})|\psi|^{\frac{p}{p-1}}\dx+\frac{\sigma_5(q-1)}{p}\int_{\Omega_\alpha}|\nabla u|^p|\psi|^p\dx+$$
$$+\frac{\sigma_5^{-\frac{q-1}{p-q+1}}C_{p,q}(p-q+1)}{p}\int_{\Omega_\alpha}(|\nabla\varphi_1|^\frac{p}{p-q+1}+|\nabla\varphi_2|^\frac{p}{p-q+1})|\psi|^\frac{p}{p-q+1}\dx+\sigma_6\int_{\Omega_\alpha}|\nabla u|^p|\psi|^p\dx+$$
$$+\overline{C}_{k,p,q}(p-1)\sigma_6^{-\frac{1}{p-1}}\int_{\Omega_\alpha}|\nabla \psi|^{\frac{p}{p-1}}(|u-\bar{u}_{B_{r}}|^{\frac{p}{p-1}}+|\varphi_1-\overline{\varphi_1}_{B_{r}}|^{\frac{p}{p-1}}+|\varphi_2-\overline{\varphi_2}_{B_{r}}|^{\frac{p}{p-1}})\dx+$$ $$+\sigma_7(q-1)\int_{\Omega_\alpha}|\nabla u|^p|\psi|^p\dx+$$
$$+\sigma_7^{-\frac{q-1}{p-q+1}}\overline{C}_{p,q}(p-q+1)\int_{\Omega_\alpha}|\nabla \psi|^{\frac{p}{p-q+1}}(|u-\bar{u}_{B_{r}}|^{\frac{p}{p-q+1}}+|\varphi_1-\overline{\varphi_1}_{B_{r}}|^{\frac{p}{p-q+1}}+|\varphi_2-\overline{\varphi_2}_{B_{r}}|^{\frac{p}{p-q+1}})\dx+$$
\begin{equation}\label{a1}
+||f\psi||_{p'}|||\psi(u-\bar{u}_{B_{r}})|+|\psi(\varphi_1-\overline{\varphi_1}_{B_{r}})|+|\psi(\varphi_2-\overline{\varphi_2}_{B_{r}})|||_{p},
\end{equation}
where for the last term we used H\"{o}lder's inequality.\\
Now, thanks to Poincaré's inequality in the last term of (\ref{a1}), we have:
$$||f\psi||_{p'}|||\psi(u-\bar{u}_{B_{r}})|+|\psi(\varphi_1-\overline{\varphi_1}_{B_{r}})|+|\psi(\varphi_2-\overline{\varphi_2}_{B_{r}})|||_{p}\leq$$
$$\leq||f\psi||_{p'}||\psi(u-\bar{u}_{B_{r}})||_{p}+||f\psi||_{p'}||\psi(\varphi_1-\overline{\varphi_1}_{B_{r}})||_{p}+||f\psi||_{p'}||\psi(\varphi_2-\overline{\varphi_2}_{B_{r}})||_{p}\leq$$
$$\leq||f\psi||_{p'}d||\nabla(\psi(u-\bar{u}_{B_{r}}))||_{p}+||f\psi||_{p'}d||\nabla(\psi(\varphi_1-\overline{\varphi_1}_{B_{r}}))||_{p}+||f\psi||_{p'}d||\nabla(\psi(\varphi_2-\overline{\varphi_2}_{B_{r}}))||_{p}=$$
$$=d||f\psi||_{p'}[||\nabla\psi(u-\bar{u}_{B_{r}})+\psi\nabla u||_p+||\nabla\psi(\varphi_1-\overline{\varphi_1}_{B_{r}})+\psi\nabla\varphi_1||_p+||\nabla\psi(\varphi_2-\overline{\varphi_2}_{B_{r}})+\psi\nabla\varphi_2||_p]\leq$$
$$\leq d||f\psi||_{p'}\Big(||\nabla\psi(u-\bar{u}_{B_{r}})||_p+||\nabla\psi(\varphi_1-\overline{\varphi_1}_{B_{r}})||_p+||\nabla\psi(\varphi_2-\overline{\varphi_2}_{B_{r}})||_p+||\psi\nabla u||_p+||\psi\nabla\varphi_1||_p+||\psi\nabla\varphi_2||_p\Big).$$
Eventually, using Young's inequality, term by term, with conjugate exponents $p$ and $\frac{p}{p-1}$, we obtain:
$$||f\psi||_{p'}|||\psi(u-\bar{u}_{B_{r}})|+|\psi(\varphi_1-\overline{\varphi_1}_{B_{r}})|+|\psi(\varphi_2-\overline{\varphi_2}_{B_{r}})|||_{p}\leq\frac{\sigma_8}{p}||\psi\nabla u||^p_p+\frac{\sigma_8}{p}||\psi\nabla \varphi_1||^p_p+\frac{\sigma_8}{p}||\psi\nabla \varphi_2||^p_p+$$
\begin{equation}\label{a2}
+\frac{\sigma_8}{p}||\nabla\psi(u-\bar{u}_{B_{r}})||^p_p+\frac{\sigma_8}{p}||\nabla\psi(\varphi_1-\overline{\varphi_1}_{B_{r}})||^p_p+\frac{\sigma_8}{p}||\nabla\psi(\varphi_2-\overline{\varphi_2}_{B_{r}})||^p_p+d^{\frac{p}{p-1}}\frac{\sigma_8^{-\frac{1}{p-1}}(p-1)}{p}||f\psi||^{p'}_{p'}.
\end{equation}
Hence, putting together (\ref{a1}) and (\ref{a2}) and suppressing the term $-||\psi\nabla u||_q^q$, we have
$$\Big(1-\frac{\sigma}{p}-\sigma_1+\frac{\sigma_1}{p}-\sigma_2-\sigma_3p+\sigma_3-\frac{\sigma_4}{p}-\frac{\sigma_5q}{p}+\frac{\sigma_5}{p}-\sigma_6-\sigma_7q+\sigma_7-\frac{\sigma_8}{p}\Big)\int_{\Omega_\alpha}|\psi|^p|\nabla u|^p\dx\leq$$
$$+\frac{\overline{C}_{k,p}(p-1)\sigma^{-\frac{1}{p-1}}}{p}\int_{\Omega_\alpha}(|\nabla\varphi_1|^{\frac{p}{p-1}}+|\nabla\varphi_2|^{\frac{p}{p-1}})|\psi|^{\frac{p}{p-1}}\dx+\frac{\sigma_1^{-(p-1)}}{p}\overline{C}_p\int_{\Omega_\alpha}(|\nabla\varphi_1|^p+|\nabla\varphi_2|^p)|\psi|^p\dx+$$
$$+\overline{C}_{k,p}(p-1)\sigma_2^{-\frac{1}{p-1}}\int_{\Omega_\alpha}|\nabla \psi|^{\frac{p}{p-1}}(|u-\bar{u}_{B_{r}}|^{\frac{p}{p-1}}+|\varphi_1-\overline{\varphi_1}_{B_{r}}|^{\frac{p}{p-1}}+|\varphi_2-\overline{\varphi_2}_{B_{r}}|^{\frac{p}{p-1}})\dx+$$
$$+\overline{C}_{p}\sigma_3^{-(p-1)}\int_{\Omega_\alpha}|\nabla \psi|^p(|u-\bar{u}_{B_{r}}|^p+|\varphi_1-\overline{\varphi_1}_{B_{r}}|^p+|\varphi_2-\overline{\varphi_2}_{B_{r}}|^p)\dx+$$
$$+\frac{C_{k,p,q}(p-1)\sigma_4^{-\frac{1}{p-1}}}{p}\int_{\Omega_\alpha}(|\nabla\varphi_1|^{\frac{p}{p-1}}+|\nabla\varphi_2|^{\frac{p}{p-1}})|\psi|^{\frac{p}{p-1}}\dx+$$
$$+\frac{\sigma_5^{-\frac{q-1}{p-q+1}}C_{p,q}(p-q+1)}{p}\int_{\Omega_\alpha}(|\nabla\varphi_1|^\frac{p}{p-q+1}+|\nabla\varphi_2|^\frac{p}{p-q+1})|\psi|^\frac{p}{p-q+1}\dx+$$
$$+\overline{C}_{k,p,q}(p-1)\sigma_6^{-\frac{1}{p-1}}\int_{\Omega_\alpha}|\nabla \psi|^{\frac{p}{p-1}}(|u-\bar{u}_{B_{r}}|^{\frac{p}{p-1}}+|\varphi_1-\overline{\varphi_1}_{B_{r}}|^{\frac{p}{p-1}}+|\varphi_2-\overline{\varphi_2}_{B_{r}}|^{\frac{p}{p-1}})\dx+$$
$$+\sigma_7^{-\frac{q-1}{p-q+1}}\overline{C}_{p,q}(p-q+1)\int_{\Omega_\alpha}|\nabla \psi|^{\frac{p}{p-q+1}}(|u-\bar{u}_{B_{r}}|^{\frac{p}{p-q+1}}+|\varphi_1-\overline{\varphi_1}_{B_{r}}|^{\frac{p}{p-q+1}}+|\varphi_2-\overline{\varphi_2}_{B_{r}}|^{\frac{p}{p-q+1}})\dx+$$
$$+\frac{\sigma_8}{p}\int_{\Omega_\alpha}|\psi|^p|\nabla \varphi_1|^p\dx+\frac{\sigma_8}{p}\int_{\Omega_\alpha}|\psi|^p|\nabla \varphi_2|^p\dx+d^{\frac{p}{p-1}}\frac{\sigma_8^{-\frac{1}{p-1}}(p-1)}{p}\int_{\Omega_\alpha}|f|^{\frac{p}{p-1}}|\psi|^{\frac{p}{p-1}}\dx+$$
$$+\frac{\sigma_8}{p}\int_{\Omega_\alpha}|\nabla\psi|^p|u-\bar{u}_{B_{r}}|^p\dx+\frac{\sigma_8}{p}\int_{\Omega_\alpha}|\nabla\psi|^p|\varphi_1-\overline{\varphi_1}_{B_{r}}|^p\dx+\frac{\sigma_8}{p}\int_{\Omega_\alpha}|\nabla\psi|^p|\varphi_2-\overline{\varphi_2}_{B_{r}}|^p\dx.$$
Lastly, choosing $\sigma=\sigma_1=\sigma_7=\frac{1}{5pq}$, $\sigma_4=\sigma_5=\frac{1}{5q}$, $\sigma_2=\sigma_3=\frac{1}{5p^2}$, $\sigma_6=\frac{1}{5p}$, $\sigma_8=\frac{1}{5}$  and considering as $\overline{C}=\overline{C}(k,p,q,\Omega)$ the maximum among all the coefficients of the integrals at the right-hand side of the last inequality, we get:
$$\int_{\Omega_\alpha}|\psi|^p|\nabla u|\dx\leq C\Big[\int_{\Omega_\alpha}|\psi|^p(|\nabla\varphi_1|^p+|\nabla\varphi_2|^p)\dx+$$ $$+\int_{\Omega_\alpha}|\psi|^{\frac{p}{p-1}}\Big(|\nabla\varphi_1|^{\frac{p}{p-1}}+|\nabla\varphi_2|^{\frac{p}{p-1}}\Big)\dx+\int_{\Omega_\alpha}|\psi|^{\frac{p}{p-q+1}}\Big(|\nabla\varphi_1|^{\frac{p}{p-q+1}}+|\nabla\varphi_2|^{\frac{p}{p-q+1}}\Big)\dx+$$
$$+\int_{\Omega_\alpha}|\nabla\psi|^p(|u-\bar{u}_{B_{r}}|^p+|\varphi_1-\overline{\varphi_1}_{B_{r}}|^p+|\varphi_2-\overline{\varphi_2}_{B_{r}}|^p)\dx+$$ $$+\int_{\Omega_\alpha}|\nabla\psi|^{\frac{p}{p-1}}(|u-\bar{u}_{B_{r}}|^{\frac{p}{p-1}}+|\varphi_1-\overline{\varphi_1}_{B_{r}}|^{\frac{p}{p-1}}+|\varphi_2-\overline{\varphi_2}_{B_{r}}|^{\frac{p}{p-1}})\dx+$$
$$+\int_{\Omega_\alpha}|\nabla\psi|^{\frac{p}{p-q+1}}(|u-\bar{u}_{B_{r}}|^{\frac{p}{p-q+1}}+|\varphi_1-\overline{\varphi_1}_{B_{r}}|^{\frac{p}{p-q+1}}+|\varphi_2-\overline{\varphi_2}_{B_{r}}|^{\frac{p}{p-q+1}})\dx+\int_{\Omega_\alpha}|f|^{\frac{p}{p-1}}|\psi|^{\frac{p}{p-1}}\dx\Big],$$
where $C=C(k,p,q,\Omega)=2\frac{p}{p-1}\overline{C}$. So, the desired relation is proved.\\

For completeness, let us show (\ref{Ineq2}). The proof is completely analogous.\\
Let us consider 
$$w:=(\varphi_1\vee g)\land\varphi_2=
\begin{cases}
\varphi_1,\text{ if } g<\varphi_1\\
g,\text{ if } \varphi_1\leq g\leq\varphi_2\\
\varphi_2,\text{ if } g>\varphi_2
\end{cases}$$ 
and 
$$v=u-\psi^p(u-w)=(1-\psi^p)u+\psi^pw.$$ 
We have that $v\in W_g^{1,p}(\Omega_\alpha)$.\\
Moreover, since $\varphi_1\leq w\leq\varphi_2$ by definition, we get that $\varphi_1\leq v\leq\varphi_2$. Then $v\in\mathcal{H}_{p}$.\\
With this choice of $v$ we have 
$$\nabla(v-u)=\nabla(-\psi^pu+\psi^pw)=p\psi^{p-1}\nabla\psi (-u+w)-\psi^p\nabla u+\psi^p\nabla w.$$ 
Analogously to the proof of part \text{(i)}, it hold that $|\nabla w|\leq |\nabla \varphi_1|+|\nabla g|+|\nabla\varphi_2|$.
Then, substituting in (\ref{Ppq}), we have:
$$0\leq\int_{\Omega_\alpha}(k^2+|\nabla u|^2)^{\frac{p-2}{2}}\nabla u\cdot[p\psi^{p-1}\nabla\psi (-u+w)-\psi^p\nabla u+\psi^p\nabla w]\dx+$$
$$+\int_{\Omega_\alpha}(k^2+|\nabla u|^2)^{\frac{q-2}{2}}\nabla u\cdot[p\psi^{p-1}\nabla\psi (-u+w)-\psi^p\nabla u+\psi^p\nabla w]\dx+
\int_{\Omega_\alpha}f\psi^p(u-w)\dx\leq$$

$$\leq-\int_{\Omega_\alpha}(k^2+|\nabla u|^2)^{\frac{p-2}{2}}|\nabla u|^2|\psi|^p\dx+\int_{\Omega_\alpha}(k^2+|\nabla u|^2)^{\frac{p-2}{2}}|\nabla u|(|\nabla\varphi_1|+|\nabla g|+|\nabla\varphi_2|)|\psi|^p\dx+$$
$$+p\int_{\Omega_\alpha}(k^2+|\nabla u|^2)^{\frac{p-2}{2}}|\nabla u||\nabla \psi||\psi|^{p-1}|u-w|\dx-\int_{\Omega_\alpha}(k^2+|\nabla u|^2)^{\frac{q-2}{2}}|\nabla u|^2|\psi|^p\dx+$$
$$+\int_{\Omega_\alpha}(k^2+|\nabla u|^2)^{\frac{q-2}{2}}|\nabla u|(|\nabla\varphi_1|+|\nabla g|+|\nabla\varphi_2|)|\psi|^p\dx+$$ $$+p\int_{\Omega_\alpha}(k^2+|\nabla u|^2)^{\frac{q-2}{2}}|\nabla u||\nabla \psi||\psi|^{p-1}|u-w|\dx+$$
$$+\int_{\Omega_\alpha}|f\psi||\psi|^{p-1}|u-w|\dx.$$
So, with the same arguments of before, we have:
$$0\leq-\int_{\Omega_\alpha}|\nabla u|^p|\psi|^p\dx+C_{k,p}\int_{\Omega_\alpha}|\nabla u|(|\nabla\varphi_1|+|\nabla g|+|\nabla\varphi_2|)|\psi|^2\dx+$$ $$+C_p\int_{\Omega_\alpha}|\nabla u|^{p-1}(|\nabla\varphi_1|+|\nabla g|+|\nabla\varphi_2|)|\psi|^p\dx+$$
$$+pC_{k,p}\int_{\Omega_\alpha}|\nabla u||\nabla \psi||\psi||u-w|\dx+pC_p\int_{\Omega_\alpha}|\nabla u|^{p-1}|\nabla \psi||\psi|^{p-1}|u-w|\dx+$$
$$-\int_{\Omega_\alpha}|\nabla u|^q|\psi|^p\dx+C_{k,q}\int_{\Omega_\alpha}|\nabla u|(|\nabla\varphi_1|+|\nabla g|+|\nabla\varphi_2|)|\psi|^2\dx+$$ $$+C_q\int_{\Omega_\alpha}|\nabla u|^{q-1}(|\nabla\varphi_1|+|\nabla g|+|\nabla\varphi_2|)|\psi|^q\dx+$$
$$+pC_{k,q}\int_{\Omega_\alpha}|\nabla u||\nabla \psi||\psi||u-w|\dx+pC_q\int_{\Omega_\alpha}|\nabla u|^{q-1}|\nabla \psi||\psi|^{q-1}|u-w|\dx+\int_{\Omega_\alpha}|f\psi||\psi||u-w|\dx\leq$$

$$\leq-\int_{\Omega_\alpha}|\nabla u|^p|\psi|^p\dx+\frac{\sigma}{p}\int_{\Omega_\alpha}|\nabla u|^p|\psi|^p\dx+\frac{\overline{C}_{k,p}(p-1)\sigma^{-\frac{1}{p-1}}}{p}\int_{\Omega_\alpha}(|\nabla\varphi_1|^{\frac{p}{p-1}}+|\nabla g|^{\frac{p}{p-1}}+|\nabla\varphi_2|^{\frac{p}{p-1}})|\psi|^{\frac{p}{p-1}}\dx+$$
$$+\frac{\sigma_1(p-1)}{p}\int_{\Omega_\alpha}|\nabla u|^p|\psi|^p\dx+\frac{\sigma_1^{-(p-1)}}{p}\overline{C}_p\int_{\Omega_\alpha}(|\nabla\varphi_1|^p+|\nabla g|^p+|\nabla\varphi_2|^p)|\psi|^p\dx+\sigma_2\int_{\Omega_\alpha}|\nabla u|^p|\psi|^p\dx+$$
$$+\overline{C}_{k,p}(p-1)\sigma_2^{-\frac{1}{p-1}}\int_{\Omega_\alpha}|\nabla \psi|^{\frac{p}{p-1}}|u-w|^{\frac{p}{p-1}}\dx+\sigma_3(p-1)\int_{\Omega_\alpha}|\nabla u|^p|\psi|^p\dx+\overline{C}_{p}\sigma_3^{-(p-1)}\int_{\Omega_\alpha}|\nabla \psi|^p|u-w|^p\dx+$$
$$-\int_{\Omega_\alpha}|\nabla u|^q|\psi|^q\dx+\frac{\sigma_4}{p}\int_{\Omega_\alpha}|\nabla u|^p|\psi|^p\dx+\frac{C_{k,p,q}(p-1)\sigma_4^{-\frac{1}{p-1}}}{p}\int_{\Omega_\alpha}(|\nabla\varphi_1|^{\frac{p}{p-1}}+|\nabla g|^{\frac{p}{p-1}}+|\nabla\varphi_2|^{\frac{p}{p-1}})|\psi|^{\frac{p}{p-1}}\dx$$
$$+\frac{\sigma_5(q-1)}{p}\int_{\Omega_\alpha}|\nabla u|^p|\psi|^p\dx+\frac{\sigma_5^{-\frac{q-1}{p-q+1}}C_{p,q}(p-q+1)}{p}\int_{\Omega_\alpha}(|\nabla\varphi_1|^\frac{p}{p-q+1}+|\nabla g|^{\frac{p}{p-q+1}}+|\nabla\varphi_2|^\frac{p}{p-q+1})|\psi|^\frac{p}{p-q+1}\dx+$$
$$+\sigma_6\int_{\Omega_\alpha}|\nabla u|^p|\psi|^p\dx+\overline{C}_{k,p,q}(p-1)\sigma_6^{-\frac{1}{p-1}}\int_{\Omega_\alpha}|\nabla \psi|^{\frac{p}{p-1}}|u-w|^{\frac{p}{p-1}}\dx+\sigma_7(q-1)\int_{\Omega_\alpha}|\nabla u|^p|\psi|^p\dx+$$
\begin{equation}\label{a1bis}
+\sigma_7^{-\frac{q-1}{p-q+1}}\overline{C}_{p,q}(p-q+1)\int_{\Omega_\alpha}|\nabla \psi|^{\frac{p}{p-q+1}}|u-w|^{\frac{p}{p-q+1}}\dx+||f\psi||_{p'}||\psi(u-w)||_{p}.
\end{equation}
Now, using Poincaré's inequality in the last term of (\ref{a1bis}), we have:
$$||f\psi||_{p'}||\psi(u-w)||_{p}\leq||f\psi||_{p'}d||\nabla\psi(u-w)+\psi\nabla u-\psi\nabla w||_{p}\leq$$
$$\leq d||f\psi||_{p'}\Big(||\nabla\psi(u-w)||_p+||\psi\nabla u||_p+||\psi\nabla w||_{p}\Big).$$
Applying Young's inequality, term by term, with conjugate exponents $p$ and $\frac{p}{p-1}$, we obtain:
\begin{equation}\label{a2bis}
||f\psi||_{p'}||\psi(u-w)||_{p}\leq\frac{\sigma_8}{p}||\psi\nabla u||^p_p+\frac{\sigma_8}{p}||\psi\nabla w||^p_p+\frac{\sigma_8}{p}||\nabla\psi(u-w)||^p_p+\frac{d^{\frac{p}{p-1}}\sigma_8^{-\frac{1}{p-1}}(p-1)}{p}||f\psi||^{p'}_{p'}.
\end{equation}
Putting together (\ref{a1bis}) and (\ref{a2bis}), and with the same choice of the parameters $\sigma$ and $\sigma_i$, $i=1\dots8$, we get
$$\int_{\Omega_\alpha}|\psi|^p|\nabla u|\dx\leq $$ $$\lq C\Big[\int_{\Omega_\alpha}|\psi|^p(|\nabla\varphi_1|^p
+|\nabla g|^p+|\nabla\varphi_2|^p)\dx+\int_{\Omega_\alpha}|\psi|^{\frac{p}{p-1}}\Big(|\nabla\varphi_1|^{\frac{p}{p-1}}+|\nabla g|^{\frac{p}{p-1}}+|\nabla\varphi_2|^{\frac{p}{p-1}}\Big)\dx+$$
$$+\int_{\Omega_\alpha}|\psi|^{\frac{p}{p-q+1}}\Big(|\nabla\varphi_1|^{\frac{p}{p-q+1}}+|\nabla g|^{\frac{p}{p-q+1}}+|\nabla\varphi_2|^{\frac{p}{p-q+1}}\Big)\dx+\int_{\Omega_\alpha}|\nabla\psi|^p|u-w|^p\dx+\int_{\Omega_\alpha}|\nabla\psi|^{\frac{p}{p-1}}|u-w|^{\frac{p}{p-1}}\dx+$$
$$+\int_{\Omega_\alpha}|\nabla\psi|^{\frac{p}{p-q+1}}|u-w|^{\frac{p}{p-q+1}}\dx+\int_{\Omega_\alpha}|f|^{\frac{p}{p-1}}|\psi|^{\frac{p}{p-1}}\dx\Big],$$
where $C=C(k,p,q,\Omega)>0.$
\endproof

\subsection{Integrability of the gradient}\label{IotG}
Now, thanks to the technical lemmas just seen, it is possible to prove the following summability result.
\begin{thm}\label{SumRes}
Let $f\in L^{p'}(\Omega_\alpha)$, $\varphi_1,\varphi_2\in W^{1,p}(\Omega_\alpha)$, $g\in W^{1,\infty}(\Omega_\alpha)$ and $u$ solution to Problem (\ref{Ppq}). Then there exists $\delta>0$, depending only on $p$ , $q$ and the $p$-thickness constant of $ \R^2\setminus\Omega_\alpha$ such that $|\nabla u|\in L^{p+\varepsilon}(\Omega_\alpha)$, $\forall\varepsilon\in(0,\delta)$.
\end{thm}
\proof Since $\Omega_\alpha$ is bounded, we can consider a ball $B_{2a}$, , such that $\Omega_\alpha\subset\subset B_{a}$. Then, let us fix $r>0$ and let us consider a ball $B_{2r}\subset B_{2a}.$Two cases are possible:\\
(i) $B_{2r}\subset\Omega_\alpha$;\\
(ii) $B_{2r}\cap(\R^2\setminus\overline{\Omega}_\alpha)\neq\emptyset$.\\
\textbf{Case (i):} Let us consider $\psi\in C_0^{\infty}(B_{2r})$, with $0\leq\psi\leq1$, $|\nabla\psi|\leq \frac{4}{r}$, and $\psi=1$ in $B_r$.
Then, by relation (\ref{Ineq1}) of Lemma \ref{Ineq} we have:
$$\int_{B_r}|\nabla u|\dx\leq C\Big[\int_{B_{2r}}\Big
(|\nabla\varphi_1|^p+|\nabla\varphi_2|^p+|\nabla\varphi_1|^{\frac{p}{p-1}}+|\nabla\varphi_2|^{\frac{p}{p-1}}+|\nabla\varphi_1|^{\frac{p}{p-q+1}}+|\nabla\varphi_2|^{\frac{p}{p-q+1}}\Big)\dx+$$
$$+2^{2p}r^{-p}\int_{B_{2r}}(|u-\bar{u}_{B_{r}}|^p+|\varphi_1-\overline{\varphi_1}_{B_{r}}|^p+|\varphi_2-\overline{\varphi_2}_{B_{r}}|^p)\dx+2^{\frac{2p}{p-1}}r^{-\frac{p}{p-1}}\int_{B_{2r}}\Big(|u-\bar{u}_{B_{r}}|^{\frac{p}{p-1}}+|\varphi_1-\overline{\varphi_1}_{B_{r}}|^{\frac{p}{p-1}}+|\varphi_2-\overline{\varphi_2}_{B_{r}}|^{\frac{p}{p-1}}\Big)\dx+$$
\begin{equation}\label{ai}
+2^{\frac{2p}{p-q+1}}r^{-\frac{p}{p-q+1}}\int_{B_{2r}}\Big(|u-\bar{u}_{B_{r}}|^{\frac{p}{p-q+1}}+|\varphi_1-\overline{\varphi_1}_{B_{r}}|^{\frac{p}{p-q+1}}+|\varphi_2-\overline{\varphi_2}_{B_{r}}|^{\frac{p}{p-q+1}}\Big)\dx+\int_{B_{2r}}|f|^{\frac{p}{p-1}}\dx\Big].
\end{equation}
Hence, considering $p=\frac{2t}{2-t}$, that is $t=\frac{2p}{p+2}$, 
by using relation (\ref{Sob2}) of Lemma \ref{Sob}, H\"{o}lder' inequality and Poincaré's inequality, by previous relation (\ref{ai}), we have:
$$\fint_{B_r}|\nabla u|^p\dx\leq C\Big[\fint_{B_{2r}}\Big
(|\nabla\varphi_1|^p+|\nabla\varphi_2|^p+|\nabla\varphi_1|^{\frac{p}{p-1}}+|\nabla\varphi_2|^{\frac{p}{p-1}}+|\nabla\varphi_1|^{\frac{p}{p-q+1}}+|\nabla\varphi_2|^{\frac{p}{p-q+1}}\Big)\dx+$$
$$+2^{2p}r^{-p}\fint_{B_{2r}}(|u-\bar{u}_{B_{r}}|^p+|\varphi_1-\overline{\varphi_1}_{B_{r}}|^p+|\varphi_2-\overline{\varphi_2}_{B_{r}}|^p)\dx+2^{\frac{2p}{p-1}}r^{-\frac{p}{p-1}}\fint_{B_{2r}}\Big(|u-\bar{u}_{B_{r}}|^{\frac{p}{p-1}}+|\varphi_1-\overline{\varphi_1}_{B_{r}}|^{\frac{p}{p-1}}+|\varphi_2-\overline{\varphi_2}_{B_{r}}|^{\frac{p}{p-1}}\Big)\dx+$$
$$+2^{\frac{2p}{p-q+1}}r^{-\frac{p}{p-q+1}}\fint_{B_{2r}}\Big(|u-\bar{u}_{B_{r}}|^{\frac{p}{p-q+1}}+|\varphi_1-\overline{\varphi_1}_{B_{r}}|^{\frac{p}{p-q+1}}+|\varphi_2-\overline{\varphi_2}_{B_{r}}|^{\frac{p}{p-q+1}}\Big)\dx+\fint_{B_{2r}}|f|^{\frac{p}{p-1}}\dx\Big]\leq$$

$$\leq C\fint_{B_{2r}}\Big
(|\nabla\varphi_1|^p+|\nabla\varphi_2|^p+|\nabla\varphi_1|^{\frac{p}{p-1}}+|\nabla\varphi_2|^{\frac{p}{p-1}}+|\nabla\varphi_1|^{\frac{p}{p-q+1}}+|\nabla\varphi_2|^{\frac{p}{p-q+1}}\Big)\dx+$$
$$+C2^{2p}r^{-p}r^p\Big(\fint_{B_{2r}}|\nabla u|^{t}\dx\Big)^{\frac{p}{t}}+C2^{\frac{2p}{p-1}}\Big(r^{-p}\fint_{B_{2r}}|u-\bar{u}_{B_{r}}|^p\dx\Big)^{\frac{1}{p-1}}+C2^{\frac{2p}{p-q+1}}\Big(r^{-p}\fint_{B_{2r}}|u-\bar{u}_{B_{r}}|^p\dx\Big)^{\frac{1}{p-q+1}}+$$
$$+C2^{2p}r^{-p}r^p\fint_{B_{2r}}|\nabla\varphi_1|^p\dx+C2^{2p}r^{-p}r^p\fint_{B_{2r}}|\nabla\varphi_2|^p\dx+C2^{\frac{2p}{p-1}}r^{-\frac{p}{p-1}}r^{\frac{p}{p-1}}\fint_{B_{2r}}|\nabla\varphi_1|^{\frac{p}{p-1}}\dx+C2^{\frac{2p}{p-1}}r^{-\frac{p}{p-1}}r^{\frac{p}{p-1}}\fint_{B_{2r}}|\nabla\varphi_2|^{\frac{p}{p-1}}\dx+$$
$$+C2^{\frac{2p}{p-q+1}}r^{-\frac{p}{p-q+1}}r^{\frac{p}{p-q+1}}\fint_{B_{2r}}|\nabla\varphi_1|^{\frac{p}{p-q+1}}\dx+C2^{\frac{2p}{p-q+1}}r^{-\frac{p}{p-q+1}}r^{\frac{p}{p-q+1}}\fint_{B_{2r}}|\nabla\varphi_2|^{\frac{p}{p-q+1}}\dx+C+\fint_{B_{2r}}|f|^{\frac{p}{p-1}}\dx\leq$$

$$\leq C\Big(\fint_{B_{2r}}|\nabla u|^{p(1-\varepsilon_1)}\dx\Big)^{\frac{1}{1-\varepsilon_1}}+C\Big[\Big(\fint_{B_{2r}}|\nabla u|^{p(1-\varepsilon_1)}\dx\Big)^{\frac{1}{1-\varepsilon_1}}\Big]^{\frac{1}{p-1}}+C\Big[\Big(\fint_{B_{2r}}|\nabla u|^{p(1-\varepsilon_1)}\dx\Big)^{\frac{1}{1-\varepsilon_1}}\Big]^{\frac{1}{p-q+1}}+$$
$$+ C\fint_{B_{2r}}\Big
(|\nabla\varphi_1|^p+|\nabla\varphi_2|^p+|\nabla\varphi_1|^{\frac{p}{p-1}}+|\nabla\varphi_2|^{\frac{p}{p-1}}+|\nabla\varphi_1|^{\frac{p}{p-q+1}}+|\nabla\varphi_2|^{\frac{p}{p-q+1}}+|f|^{\frac{p}{p-1}}\Big)\dx$$
with $0<\varepsilon_1\leq\frac{p}{p+2}$.\\
Since, by Young's inequality, we have that 
\begin{equation}\label{young}
a^{\frac{1}{p}}\leq \frac{a}{p^{\frac{1}{p}}}+\frac{p-1}{p}\leq a+1,\,\,\forall a\geq0 \text{ and } \forall p\geq1,
\end{equation}
we get
\begin{equation}\label{ao}
\fint_{B_r}|\nabla u|^p\dx\leq C\Big(\fint_{B_{2r}}|\nabla u|^{p(1-\varepsilon_1)}\dx\Big)^{\frac{1}{1-\varepsilon_1}}+\fint_{B_{2r}}\Big[C^{\frac{1}{p}}\Big
(|\nabla\varphi_1|+|\nabla\varphi_2|+|\nabla\varphi_1|^{\frac{1}{p-1}}+|\nabla\varphi_2|^{\frac{1}{p-1}}+|\nabla\varphi_1|^{\frac{1}{p-q+1}}+|\nabla\varphi_2|^{\frac{1}{p-q+1}}+|f|^{\frac{1}{p-1}}+2^{\frac{1}{p}}\Big)\Big]^p\dx.
\end{equation}

\textbf{Case (ii):}
First of all, let $Ext_J(g)$ be the extension of $g$ on $\R^2$, which is possible since $\Omega_\alpha$ is a $(\epsilon,\delta)$-domain (see \cite{C} and \cite{J}, f.i.). So, we extend $u$, $w$ and $f$ on $\R^2$ setting that they are equal to $Ext_J(g)$ in $\R^2\setminus\Omega_\alpha$.\\ 
Let us consider $\psi$ as in case \textit{(i)} and let $D=B_{2r}\cap\Omega_\alpha$ be. By relation (\ref{Ineq2}) of Lemma \ref{Ineq}, we have:
$$\int_{D}|\psi|^p|\nabla u|\dx\leq C\Big[\int_{D}|\psi|^p(|\nabla\varphi_1|^p
+|\nabla g|^p+|\nabla\varphi_2|^p)\dx+\int_{D}|\psi|^{\frac{p}{p-1}}\Big(|\nabla\varphi_1|^{\frac{p}{p-1}}+|\nabla g|^{\frac{p}{p-1}}+|\nabla\varphi_2|^{\frac{p}{p-1}}\Big)\dx+$$
$$+\int_{D}|\psi|^{\frac{p}{p-q+1}}\Big(|\nabla\varphi_1|^{\frac{p}{p-q+1}}+|\nabla g|^{\frac{p}{p-q+1}}+|\nabla\varphi_2|^{\frac{p}{p-q+1}}\Big)\dx+\int_{B_{2r}}|\nabla\psi|^p|u-w|^p\dx+\int_{B_{2r}}|\nabla\psi|^{\frac{p}{p-1}}|u-w|^{\frac{p}{p-1}}\dx$$
\begin{equation}\label{aii}
+\int_{B_{2r}}|\nabla\psi|^{\frac{p}{p-q+1}}|u-w|^{\frac{p}{p-q+1}}\dx+\int_{D}|f|^{\frac{p}{p-1}}|\psi|^{\frac{p}{p-1}}\dx\Big].
\end{equation}
Let us consider $q=p(1-\varepsilon_2)$, with 
$0<\varepsilon_2\leq\min\Big\{\frac{p-2}{2},\frac{1}{2}\Big\}$. So, defining
$$\alpha=\begin{cases}
\frac{2}{2-q}, \text{ if } 1<q<2\\
2, \text{ if } q\geq2\,\,\,,
\end{cases}$$ 
we have $\alpha q\geq p$ and by lemma \ref{Cap} we get:
$$\Big(\fint_{B_{2r}}|u-w|^p|\nabla\psi|^p\dx\Big)^{\frac{1}{p}}\leq4r^{-1}\Big(\fint_{B_{2r}}|u-w|^p\dx\Big)^{\frac{1}{p}}\leq$$
\begin{equation}\label{aiii}
\leq 4r^{-1}\Big(\fint_{B_{2r}}|u-w|^{\alpha q}\dx\Big)^{\frac{1}{\alpha q}}\leq C\Big(\frac{r^{2-q}}{\capa_q(N(u-g);B_{4r})}\fint_{B_{2r}}|\nabla(u-w)|^{q}\dx\Big)^{\frac{1}{q}},
\end{equation}
with $N(u-w)=\{x\in B_{2r}\,:\,u=w\}$; in $\R^2\setminus\overline{\Omega}_\alpha$ we have $u=w=Ext(g)$. We point out that $w=g$ on $\partial\Omega_\alpha$, since $\varphi_1\leq g\leq\varphi_2$ on $\partial\Omega_\alpha$.\\
Since $p>2\Longrightarrow q>1$, so by Remark \ref{pthi}, we have that
\begin{equation}\label{aiv}
{\capa}_q(N(u-w);B_{4r})\geq{\capa}_q(B_{2r}\setminus\Omega_\alpha;B_{4r})\geq Cr^{2-q}.
\end{equation}
Since, by H\"{o}lder's inequality, we have that
\begin{equation}\label{av}
\fint_{B_{2r}}|u-w|^\frac{p}{p-1}|\nabla\psi|^\frac{p}{p-1}\dx\leq\Big(\fint_{B_{2r}}|u-w|^p|\nabla\psi|^p\dx\Big)^{\frac{1}{p-1}};\fint_{B_{2r}}|u-w|^\frac{p}{p-q+1}|\nabla\psi|^\frac{p}{p-q+1}\dx\leq\Big(\fint_{B_{2r}}|u-w|^p|\nabla\psi|^p\dx\Big)^{\frac{1}{p-q+1}},
\end{equation}
then, by relations (\ref{young}), (\ref{aii}), (\ref{aiii}), (\ref{aiv}) and (\ref{av}), we get:
$$r^{-2}\int_{D}|\psi|^p|\nabla u|\dx\leq $$
$$\leq C\Big\{r^{-2}\int_{D}|\psi|^p(|\nabla\varphi_1|^p
+|\nabla g|^p+|\nabla\varphi_2|^p)\dx+r^{-2}\int_{D}|\psi|^{\frac{p}{p-1}}\Big(|\nabla\varphi_1|^{\frac{p}{p-1}}+|\nabla g|^{\frac{p}{p-1}}+|\nabla\varphi_2|^{\frac{p}{p-1}}\Big)\dx+$$
$$+r^{-2}\int_{D}|\psi|^{\frac{p}{p-q+1}}\Big(|\nabla\varphi_1|^{\frac{p}{p-q+1}}+|\nabla g|^{\frac{p}{p-q+1}}+|\nabla\varphi_2|^{\frac{p}{p-q+1}}\Big)\dx+\Big(\fint_{B_{2r}}|\nabla(u-w)|^{q}\dx\Big)^{\frac{p}{q}}+\Big[\Big(\fint_{B_{2r}}|\nabla(u-w)|^{q}\dx\Big)^{\frac{p}{q}}\Big]^{\frac{1}{p-1}}$$
$$+\Big[\Big(\fint_{B_{2r}}|\nabla(u-w)|^{q}\dx\Big)^{\frac{p}{q}}\Big]^{\frac{1}{p-q+1}}+r^{-2}\int_{D}|f|^{\frac{p}{p-1}}|\psi|^{\frac{p}{p-1}}\dx\Big\}\leq$$
$$\leq C\Big\{r^{-2}\int_{D}|\psi|^p(|\nabla\varphi_1|^p
+|\nabla g|^p+|\nabla\varphi_2|^p)\dx+r^{-2}\int_{D}|\psi|^{\frac{p}{p-1}}\Big(|\nabla\varphi_1|^{\frac{p}{p-1}}+|\nabla g|^{\frac{p}{p-1}}+|\nabla\varphi_2|^{\frac{p}{p-1}}\Big)\dx+$$
$$+r^{-2}\int_{D}|\psi|^{\frac{p}{p-q+1}}\Big(|\nabla\varphi_1|^{\frac{p}{p-q+1}}+|\nabla g|^{\frac{p}{p-q+1}}+|\nabla\varphi_2|^{\frac{p}{p-q+1}}\Big)\dx+3\Big(\fint_{B_{2r}}|\nabla(u-w)|^{q}\dx\Big)^{\frac{p}{q}}+2+r^{-2}\int_{D}|f|^{\frac{p}{p-1}}|\psi|^{\frac{p}{p-1}}\dx\Big\}.$$
Observing that $|D|\leq|B_{2r}|=4\pi r^2$ and
$$\Big(\fint_{B_{2r}}|\nabla(u-w)|^{q}\dx\Big)^{\frac{p}{q}}\dx= C\Big(r^{-2}\int_{D}|\nabla(u-w)|^{q}\dx\Big)^{\frac{p}{q}}\leq C\Big(r^{-2}\int_{D}|\nabla u|^{q}\dx\Big)^{\frac{p}{q}}+C\Big(r^{-2}\int_{D}|\nabla w|^{q}\dx\Big)^{\frac{p}{q}}\leq$$
$$ \leq C\Big(r^{-2}\int_{D}|\nabla u|^{q}\dx\Big)^{\frac{p}{q}}+Cr^{-2}\int_{D}|\nabla w|^{p}\dx\leq C\Big(r^{-2}\int_{D}|\nabla u|^{q}\dx\Big)^{\frac{p}{q}}+r^{-2}C\int_{D}(|\nabla\varphi_1|^{p}+|\nabla g|^{p}+|\nabla\varphi_2|^{p})\dx,$$
finally, we get
$$r^{-2}\int_{B_r\cap\Omega_\alpha}|\nabla u|^p\dx\leq C\Big(r^{-2}\int_{D}|\nabla u|^{p(1-\varepsilon_2)}\dx\Big)^{\frac{1}{1-\varepsilon_2}}+$$
\begin{equation}\label{avi}
+r^{-2}\int_{D}\Big[C^{\frac{1}{p}}\Big
(|\nabla\varphi_1|+|\nabla\varphi_2|+|\nabla\varphi_1|^{\frac{1}{p-1}}+|\nabla\varphi_2|^{\frac{1}{p-1}}+|\nabla\varphi_1|^{\frac{1}{p-q+1}}+|\nabla\varphi_2|^{\frac{1}{p-q+1}}+|f|^{\frac{1}{p-1}}+2^{\frac{1}{p}}\Big)\Big]^p\dx.
\end{equation}
Now, let 
$$g(x)=\begin{cases}
|\nabla u|^{p(1-\varepsilon)}, & x\in\Omega_\alpha\\
0, & x\in\R^2\setminus\Omega_{\alpha}\,\,\, ,
\end{cases}$$
$$h(x)=$$$$\begin{cases}
\Big[C^{\frac{1}{p}}\Big
(|\nabla\varphi_1|+|\nabla\varphi_2|+|\nabla\varphi_1|^{\frac{1}{p-1}}+|\nabla\varphi_2|^{\frac{1}{p-1}}+|\nabla\varphi_1|^{\frac{1}{p-q+1}}+|\nabla\varphi_2|^{\frac{1}{p-q+1}}+|f|^{\frac{1}{p-1}}+2^{\frac{1}{p}}\Big)\Big]^{p(1-\varepsilon)}, & x\in\Omega_\alpha\\
0, & x\in\R^2\setminus\Omega_{\alpha}
\end{cases}$$
and $s=\frac{1}{1-\varepsilon}$, with $\varepsilon=\min\{\varepsilon_1,\varepsilon_2\}$ such that (\ref{ao}) and (\ref{avi}) hold. 
Thus, we get the following
\begin{equation}\label{avii}
\fint_{B_r}g^s\leq C\Big(\fint_{B_{2r}}g\dx\Big)^{s}+\fint_{B_{2r}}h^{s}dx,
\end{equation}
for any $B_{2r}\subset B_{2a}.$
Hence, by lemma \ref{Som} we have the thesis and the proof is over.
\endproof

\section{Asymptotics as $p\to\infty$ and $n$ fixed}\label{ABU}
\setcounter{section}{5} 
In this section we focus our attention on the asymptotic behavior of the solutions. In particular, we prove an analogous result to the one presented in \cite{BJ} for the homogeneous case without obstacles (see \cite{BDM}, \cite{CF}, \cite{F}, \cite{IL} and \cite{MRT}  for the $p$-Laplacian case).

\begin{thm}\label{AsymPtoI}
Let us assume $f\in L^1(\Omega_\alpha)$, $\varphi_1,\varphi_2\in C(\overline{\Omega}_\alpha)$, $g\in W^{1,\infty}(\Omega_\alpha)$ (with Lipschitz constant $L$),
$$\mathcal{H}=\Big\{v\in W_g^{1,\infty}(\Omega_\alpha):\, \varphi_1\leq v\leq\varphi_2\text{ in }\Omega_\alpha,\,||\nabla v||_\infty\leq\max\{1,\sqrt{L^2+k^2}\}\Big\}\neq\emptyset$$ 
and $u_{p,q}$ the solution to Problem (\ref{Ppq}).\\
Then for any subsequence
$u_{p_k,q}$ there exists a subsubsequence, still denoted with $u_{p_k,q}$, such that, as $k\to\infty$, $u_{p_k,q}\to u_{\infty,q}$ uniformly in $C(\overline{\Omega}_\alpha)$ and weakly in
$W^{1,t}(\Omega_\alpha)$
 where the limit $ u_{\infty,q}$  belongs to $W^{1,t}(\Omega_\alpha)$ and
verifies   $$||\nabla u_{\infty,q}||_\infty\leq\max\{1,\sqrt{L^2+k^2}\}.$$ 
Moreover,

 if $L^2+k^2\leq1$, then $u_{\infty,q}$ is the unique solution to the following variational problem
\begin{equation*}\label{Pq1} 
\min_{v\in\mathcal{H}}J_{q}(v), \text{ with } J_{q}(v)=\frac{1}{q}\int_{\Omega_\alpha}(k^2+|\nabla v|^2)^\frac{q}{2}\dx-\int_{\Omega_\alpha}fv\dx\,\,\tag{ $\mathcal{P}_{q}$}
\end{equation*} and 
if $L^2+k^2>1$ then  $u_{\infty,q}$  is a minimal Lipschitz extension, that is, $u_{\infty,q}$  is a solution to
\begin{equation*}\label{PqL1} 
\min_{v\in\mathcal{H}}||\nabla v||_\infty. \,\,\tag{ $\mathcal{P}_{q.L}$}
\end{equation*}
\end{thm}
\proof
Let $w\in\mathcal{H}$ and $u_{p,q}$ solution to Problem (\ref{Ppq}). By the equivalence between this problem and Problem (\ref{PMpq}), we have
$$\frac{1}{p}||\nabla u_{p,q}||^p_p\leq \frac{1}{p}\int_{\Omega_\alpha}(k^2+|\nabla u_{p,q}|^2)^\frac{p}{2}\dx+\frac{1}{q}\int_{\Omega_\alpha}(k^2+|\nabla u_{p,q}|^2)^\frac{q}{2}\dx\leq $$
$$\frac{1}{p}\int_{\Omega_\alpha}(k^2+|\nabla v|^2)^\frac{p}{2}\dx+\frac{1}{q}\int_{\Omega_\alpha}(k^2+|\nabla v|^2)^\frac{q}{2}\dx+\int_{\Omega_\alpha}fu_{p,q}\dx-\int_{\Omega_\alpha}fv\dx\leq$$
$$\leq |\Omega_\alpha|\Big[\frac{(k^2+L^2)^\frac{p}{2}}{p}+\frac{(k^2+L^2)^\frac{q}{2}}{q}\Big]+C(f,\varphi_1,\varphi_2).$$
So
$$||\nabla u_{p,q}||_p\leq p^{\frac{1}{p}}\Big\{|\Omega_\alpha|\Big[\frac{(k^2+L^2)^\frac{p}{2}}{p}+\frac{(k^2+L^2)^\frac{q}{2}}{q}\Big]+C\Big\}^{\frac{1}{p}},$$
that is $\{u_{p,q}\}_{p>q}$ is bounded in $W^{1,p}(\Omega_\alpha)$.
Moreover, we get
$$\limsup_{p\to\infty}||\nabla u_{p,q}||_p\leq\max\{1,\sqrt{L^2+k^2}\}.$$
Now, for any $p>t>q$, by H\"{o}lder's inequality, it holds
$$||\nabla u_{p,q}||_t\leq|\Omega_\alpha|^{\frac{1}{t}-\frac{1}{p}}||\nabla u_{p,q}||_p.$$
Then, we obtain
$$\limsup_{p\to\infty}||\nabla u_{p,q}||_t\leq|\Omega_\alpha|^{\frac{1}{t}}\max\{1,\sqrt{L^2+k^2}\}.$$
Hence, by Ascoli-Arzelà compactness criterion there exist a subsequence $\{u_{p_k,q}\}_{k\in\N}$ converging to $u_{\infty,q}$ weakly in $W^{1,t}(\Omega_\alpha)$ and uniformly in $\overline{\Omega}_\alpha$.
Thus, we get
$$||\nabla u_{\infty,q}||_\infty\leq\lim_{t\to\infty}||\nabla u_{\infty,q}||_t\leq \lim_{t\to\infty}\liminf_{p\to\infty}||\nabla u_{p,q}||_t\leq\max\{1,\sqrt{L^2+k^2}\},$$ 
that is $u_{\infty,q}\in\mathcal{H}$.\\
Finally, if $k^2+L^2\leq1$, for any $v\in\mathcal{H}$, we obtain
$$\frac{1}{q}\int_{\Omega_\alpha}(k^2+|\nabla u_{p,q}|^2)^\frac{p}{2}\dx-\int_{\Omega_\alpha}fu_{p,q}\dx\leq\frac{1}{p}\int_{\Omega_\alpha}(k^2+|\nabla v|^2)^\frac{p}{2}\dx+$$
$$+\frac{1}{q}\int_{\Omega_\alpha}(k^2+|\nabla v|^2)^\frac{q}{2}\dx-\int_{\Omega_\alpha}fv\dx\leq$$
$$\leq\frac{|\Omega_\alpha|}{p}+\frac{1}{q}\int_{\Omega_\alpha}(k^2+|\nabla v|^2)^\frac{q}{2}\dx-\int_{\Omega_\alpha}fv\dx.$$
Hence, passing to the limit as $p_k\to\infty$, $p_k$ subsequence of $p$, we get that $u_{\infty,q}$ solves
$$\min_{\substack{v\in\mathcal{H}}}J_{q}(v), \text{ with } J_{q}(v)=\frac{1}{q}\int_{\Omega_\alpha}(k^2+|\nabla v|^2)^\frac{q}{2}\dx-\int_{\Omega_\alpha}fv\dx.$$
In the case $k^2+L^2>1$, by the fact that $||\nabla u_{\infty,q}||_{\infty}\leq \sqrt{k^2+L^2}$, it follows that $u_{\infty,q}$ solves
$$\min_{v\in\mathcal{H}}||\nabla v||_\infty.$$
\endproof

Until now, we have considered the problem in the setting of fractal boundary domain $\Omega_\alpha$. However, it is possible to consider the corresponding approximating problems, that is the problems on the pre-fractal approximating domains $\Omega_\alpha^n$.\\
We point out that the introduction of the following Problems (\ref{Ppqn}), besides being interesting in itself, is justified also from the possibility to perform on it numerical analysis.

Let $p>q$ be, with $q\in[2,\infty)$ fixed as before. Given $f_n\in L^{1}(\Omega_\alpha)$, $g\in W^{1,\infty}(\Omega_\alpha)$ and $\varphi_{1,n},\varphi_{2,n}\in C(\overline{\Omega}_\alpha)$, let us introduce the following problems: find  $u_{p,q,n}\in\mathcal{H}_{p,n}$ such that 
\begin{equation*}\label{Ppqn} 
a_{p,n}(u_{p,q,n},v-u_{p,q,n})+a_{q,n}(u_{p,q,n},v-u_{p,q,n})-\int_{\Omega^n_\alpha}f_n(v-u_{p,q,n})\dx \geqslant 0,\forall v\in \mathcal{H}_{p,n},\,\tag{ $\mathcal{P}_{p,q,n}$}
\end{equation*}
where
\begin{equation}
 a_{p,n}(u,v)=\int_{\Omega^n_\alpha}(k^2+|\nabla u|^2)^{\frac{p-2}{2}}\nabla u\nabla v\dx
\end{equation}
and 
\begin{equation}\label{convyn}
\mathcal{H}_{p,n}=\{v\in W_{g}^{1,p}(\Omega^n_\alpha):\,\varphi_{1,n}\leq v\leq\varphi_{2,n}\text{ in }\Omega^n_\alpha\}
\end{equation}
is assumed to be non-empty.\\
Thanks to the equivalence between this problem and the analogous of Problem (\ref{PMpq}) and Proposition \ref{UniqPpq}, adapted to $\Omega_{\alpha}^n$, we have existence and  uniqueness of the solution for the approximating Problem (\ref{Ppqn}).\\
For this problem, the following analogous result to Theorem \ref{AsymPtoI} holds.
\begin{thm}\label{AsymPtoIn}
Let us assume $f_n\in L^1(\Omega_\alpha)$, $\varphi_{1,n},\varphi_{2,n}\in C(\overline{\Omega}_\alpha)$,  $g\in W^{1,\infty}(\Omega_\alpha)$ (with Lipschitz constant $L$),
$$\mathcal{H}_n=\Big\{v\in W_g^{1,\infty}(\Omega^n_\alpha):\, \varphi_{1,n}\leq v\leq\varphi_{2,n}\text{ in }\Omega^n_\alpha,\,||\nabla v||_{\infty,\Omega^n_\alpha}\leq\max\{1,\sqrt{L^2+k^2}\}\Big\}\neq\emptyset$$ 
and $u_{p,q,n}$ the solution to Problem (\ref{Ppqn}).\\
Then for any subsequence
$u_{p_k,q,n}$ there exists a subsubsequence, still denoted with $u_{p_k,q,n}$, such that, as $k\to\infty$, $u_{p_k,q,n}\to u_{\infty,q,n}$ uniformly in $C(\overline{\Omega}^n_\alpha)$ and weakly in
$W^{1,t}(\Omega^n_\alpha)$, being $u_{\infty,q,n}$  solution of 
\begin{equation*}\label{Pqn} 
\min_{v\in\mathcal{H}_n}J_{q}(v), \text{ with } J_{q}(v)=\frac{1}{q}\int_{\Omega^n_\alpha}(k^2+|\nabla v|^2)^\frac{q}{2}\dx-\int_{\Omega^n_\alpha}fv\dx, \text{ if } k^2+L^2\leq1,\,\,\tag{ $\mathcal{P}_{q,n}$}
\end{equation*}
\begin{equation*}\label{PqLn} 
\min_{v\in\mathcal{H}_n}||\nabla v||_{\infty,\Omega^n_\alpha}, \text{ if } k^2+L^2>1\,\,\tag{ $\mathcal{P}_{q,L,n}$}
\end{equation*}
\end{thm}

\section{Asymptotics as $n\to\infty$ and $p$ fixed}

As recalled is Section  \ref{FpFB}, the sets $\Omega_\alpha^n$ give at the limit $\Omega_\alpha$. Then, it makes sense to ask whether the solutions to the approximating problems converge in some sense to a solution of the corresponding problem on $\Omega_\alpha$.\\
As far as we know, this double study on convergence, that is the analysis of the behavior with respect to $n$ as well as on $p$, was done for first in \cite{CF}, and then in \cite{F}. Nevertheless, the study of the asymptotic behavior with respect to $n$ was done by many authors for different problems (see, f.i., \cite{C}, \cite{CV3}, \cite{CV1} and \cite{MV}).\\

In order to prove the following result, let us consider $u_{p,q,n}$ solutions of Problems (\ref{Ppqn}) and define
\begin{equation}\label{Exti2}
\tilde{u}_{p,q,n}(x):=\begin{cases}
u_{p,q,n}(x), & x\in\overline{\Omega}^n_{\alpha}\\
g(x), & x\in\overline{\Omega}_{\alpha}\setminus\Omega^n_{\alpha}\,\, .
\end{cases}
\end{equation}

\begin{thm}\label{AsymRes}
Let $f_n,f\in L^{p'}(\Omega_\alpha)$, $g\in W^{1,\infty}(\Omega_\alpha)$, $\varphi_{i,n},\varphi_i\in W^{1,p}(\Omega_\alpha)$, for $i=1,2$. Moreover let us assume $\mathcal{H}_{p,n}\neq\emptyset$, $\mathcal{H}_{p}\neq\emptyset$ and, as $n\to\infty$,
$$f_n\to f \text{ in } L^1(\Omega_\alpha)\,\,\, \text{ and }\,\,\, \varphi_{i,n}\to\varphi_{i},\,\, i=1,2,\,\, \text{ in } W^{1,p}(\Omega_\alpha).$$ 
Then the sequence $\tilde{u}_{p,q,n}$ defined in (\ref{Exti2}) strongly converge, as $n\to\infty$, in $W^{1,p}(\Omega_\alpha)$ to the solution to Problem (\ref{Ppq}).
\end{thm}

Before the proof we need some preliminary results.

Let us recall how to construct a suitable array of fibers $\Sigma^n$ around the boundary of $\Omega_\alpha^n$ (see, for instance, \cite{CV1} and \cite{MV}).To show how this construction works, we start considering the open triangle of vertices  $A(0,0)$, $B(1,0)$ and $C(1/2,-\sqrt{3}/2)$.\\ Denoting with $T^+_0$  the open triangle of vertices  $A(0,0)$, $B(1,0)$ and $D^+(1/2,\delta_+/2)$, with $\delta_+= \tan(\frac{\vartheta}{2})$ and $\vartheta$ the rotation angle defined in (\ref{teta}), we have that $T^+_0$ satisfies the open set condition with respect to the family of maps $\Psi_{\alpha}$; that is $\psi_{i|n,\alpha}(T^+_0)\subset T^+_0$ for every $i|n$ and $\psi_{i|n,\alpha}(T^+_0)\cap \psi_{j|n,\alpha}(T^+_0)=\emptyset $ for every $i|n\neq j|n$. Furthermore, with $T^-_0$ we denote  the open triangle of vertices $A(0,0)$, $B(1,0)$ and $S^-(1/2,-\delta_-/2),$ where $\delta_-=\tan(\vartheta^-)$, with $0<\vartheta^-\leq\min\{\pi/2-\vartheta,\vartheta/2\}$. So, we obtain the fiber $\Sigma^0_1$ corresponding to the side $AB$ setting 
$$\Sigma_1^0=T^+_0\bigcup T^-_0\bigcup K^0.$$

Now, applying the maps $\psi_{i|n}=\psi_{i_1}\circ \psi_{i_2} \circ \cdots \circ\psi_{i_n}$, for any integer $n>0$, $\Sigma_1^0$ is iteratively transformed into increasingly fine arrays. In particular, for every $n\geq 1$, we set 
$$\Sigma^n_{1}= \Sigma^n_{1,+}\bigcup\Sigma^n_{1,-}\bigcup K^n$$ 
with
\begin{equation}\label{sigma+}
\Sigma^n_{1,+}=\bigcup_{i|n}\Sigma_{1,+}^{i|n}\,\,,\quad  \Sigma_{1,+}^{i|n}=\psi_{i|n}(T^+_0)\,\,,
\end{equation}
\begin{equation}\label{sigma-}
\Sigma^n_{1,-}=\bigcup_{i|n}\Sigma_{1,-}^{i|n}\,\,,\quad
\Sigma_{1,-}^{i|n}=\psi_{i|n}(T^-_0)\,.
\end{equation}

Denoting by $\Sigma^n_{2,+}$, $\Sigma^n_{3,+}$, $\Sigma^n_{2,-}$ and  $\Sigma^n_{3,-}$ the corresponding arrays of fibers obtained applying the same procedure to the others sides of the starting domain, we have
\begin{equation}\label{sigmas}   \Sigma^n=\bigcup_{j=1,2,3} \Sigma^n_{j},\qquad 
\Sigma_{+}^{n}= \bigcup_{j=1,2,3} \Sigma^n_{j,+},\qquad
\Sigma_{-}^n=\bigcup_{j=1,2,3} \Sigma^n_{j,-}.
\end{equation}

Hence, we define the sets
\begin{equation}\label{domExtInt}
\hat{\Omega}_\alpha^n= int(\overline{\Omega}_\alpha^n\bigcup \Sigma^n_{+})\,\,\,\text{ and }\,\,\, \breve{\Omega}_\alpha^n=\Omega_\alpha^n\setminus\overline{\Sigma}^n_{-}.
\end{equation}
In particular, we observe that for these sets it holds that
\begin{equation}\label{inclu}  
\breve{\Omega}_\alpha^n \subset \Omega_\alpha^n\subset  \hat{\Omega}_\alpha^n,
\end{equation}
\begin{equation}\label{inclu12}  
\hat{\Omega}_\alpha^{n+1} \subset \hat{\Omega}_\alpha^{n}\qquad \text{ and }\qquad  \breve{\Omega}_\alpha^n\subset \breve{\Omega}_\alpha^{n+1}. 
\end{equation}

\begin{figure}[H]
\begin{center}
	\includegraphics[height=4.0cm]{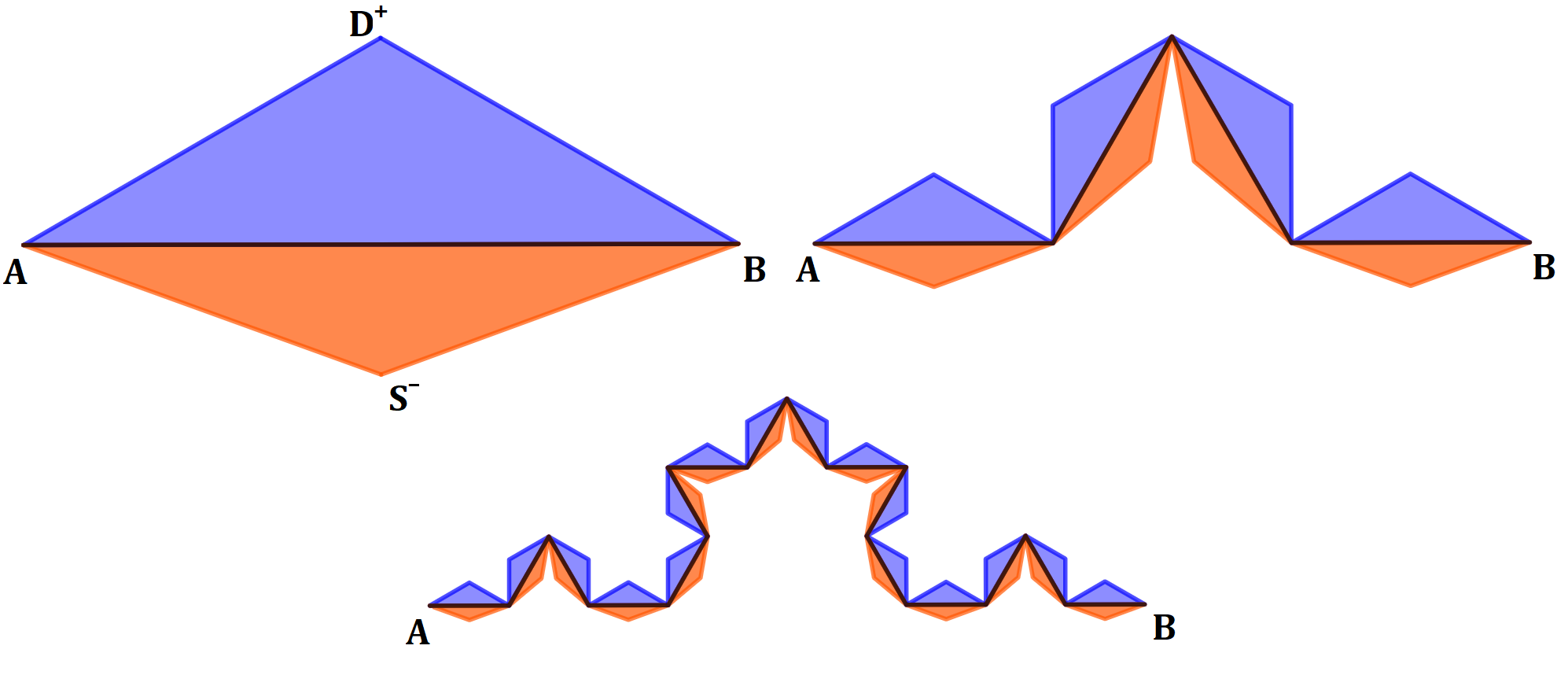}
\end{center}
\caption{$\Sigma^0_{1}$, $\Sigma^1_{1}$ and $\Sigma^2_{1}$.}
\label{Fibry}
\end{figure}

Figure \ref{Fibry} shows first iterations of the procedure just described.

\begin{figure}[H]
\begin{center}
     \includegraphics[height=5.0cm]{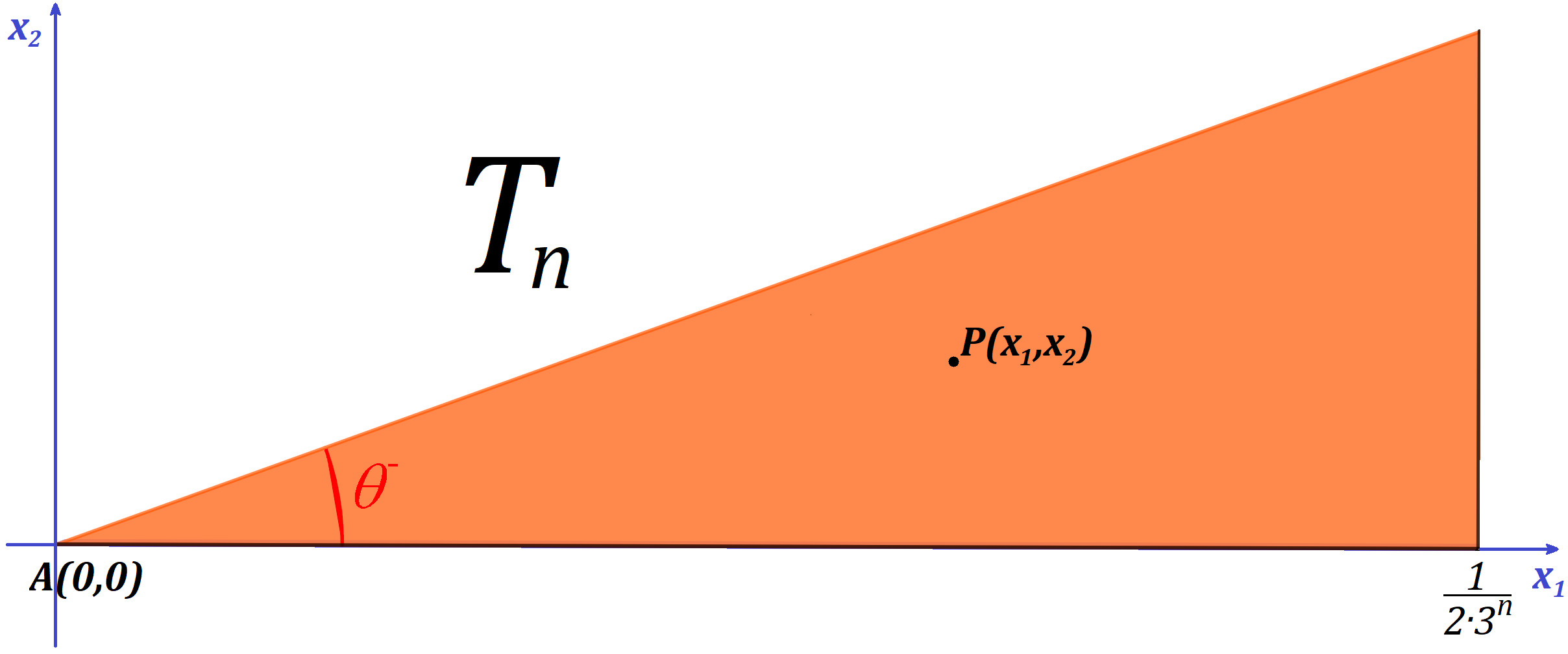}
\end{center}
\caption{The half-fiber $T_n$.}
\label{Tn}
\end{figure}

 Let us introduce a suitable function, which plays the role of coefficient of a convex combination. It allows us to construct an appropriate sequence of functions.\\
For every $n\in\N$, given $P(x_1,x_2)\in\Sigma^n_-$, we define $P_\perp(x_1^\perp,x_2^\perp)\in \partial\Omega^n_\alpha$ as the orthogonal projection of $(x_1,x_2)$ on $\partial\Omega^n_\alpha$. Then, with $P_-(x_1^-,x_2^-)$ we indicate the intersection of the straight line passing through $(x_1,x_2)$ and $(x_1^\perp,x_2^\perp)$ with $\partial\Sigma^n_-\setminus K^n$, where the symbol $-$ indicates the inner intersection. \\
Hence, we define
\begin{equation} \label{CoeLam}
\lambda_n(x)=
\begin{cases}
1,\,\,& x\in\overline{\breve{\Omega}}_{\alpha}^n\\
\frac{|x_1^{\perp}-x_1|+|x_2^{\perp}-x_2|}{|x_1^{\perp}-x_1^-|+|x_2^{\perp}-x_2^-|},\,\,& x\in\Sigma^n_-\\
0,\,\,& x\in\overline{\Omega}_\alpha\setminus\Omega_{\alpha}^n
\end{cases}
\end{equation}
with $x(x_1,x_2)$.\\
Now, let us state and prove a result which will play a central role in the analysis of the asymptotic behavior, as $n\to\infty$.
\begin{thm} \label{OpeW}
Let $u$ be in $W_0^{1,r}(\Omega_\alpha)$, $r>2$. Then, the function $w_n(x)=\lambda_n(x)\cdot u(x)$, where $\lambda_n(x)$ is defined in (\ref{CoeLam}), has the following properties:\\
\begin{eqnarray*}
	(i) & \,\,\, & w_n(x)\in W_0^{1,s}(\Omega_\alpha), \forall 2<s<r;\\
	(ii) & \,\,\, & ||w_n||_{1,s,\Omega_\alpha}\leq C, \text{ with } C \text{ independent on } n;\\
	(iii) & \,\,\, & w_n\to u \text{ in } W_0^{1,s}(\Omega_\alpha), \text{ as } n\to\infty.
\end{eqnarray*}
\end{thm}
\proof
To prove \textit{(i)} and \textit{(ii)}, let us consider $2<s<p$ and $||w_n||_{1,s,\Omega_\alpha}=||w_n||^s_{1,s}$.
$$||w_n||^s_{1,s}= ||w_n||^s_s+||\nabla w_n||^s_s\leq ||u||^s_s+||(\nabla\lambda_n)u+\lambda_n\nabla u||^s_s\leq$$
$$\leq ||u||^s_s+2^{s-1}(||(\nabla\lambda_n)u||_s^s+||\lambda_n\nabla u||^s_s)\leq ||u||^s_s+2^{s-1}||\nabla u||^s_s+2^{s-1}||(\nabla\lambda_n)u||_s^s\leq$$
$$\leq 2^{s-1}\Big(||u||^s_{1,s}+\int_{\Omega_\alpha}|(\nabla\lambda_n)u|^s\dx\Big).$$
By definition of $\lambda_n$, we get:
\begin{equation}\label{norm}
||w_n||^s_{1,s}\leq2^{s-1}\Big(||u||^s_{1,s}+\int_{\Sigma_-^n}|(\nabla\lambda_n)u|^s\dx\Big)=2^{s-1}\Big(||u||^s_{1,s}+\sum_{i=1}^{3\cdot4^n}\int_{T^n_i}|(\nabla\lambda_n)u|^s\dx\Big),
\end{equation}
were $T^n_i$ indicate the i-th triangle of the internal fiber (see the definition of $\Sigma^{i|n}_{1,-}$).\\
Now, let us focus our attention on an half-fiber triangle (that we indicate with $T_n$) having a vertex on the point $A(0,0)$ and a side on the abscissa axis (see Figure \ref{Tn}). By rotation and translation, the conclusions hold also for the other verteces of $\Omega^n_\alpha$.

In our model case $\lambda_n(x)$ and $T_n$ have the following forms:
\begin{equation}\label{lambdaTn}
\lambda_n(x)=\frac{x_2}{x_1a},\,\, x\in\Sigma^n_-,
\end{equation} 
with $a=\tan\theta^-$, and
$$T_n=\{(x_1,x_2)\in\R^2\,:\, 0\leq x_1\leq\frac{1}{2\cdot3^n},\,0\leq x_2\leq ax_1\}.$$
By (\ref{lambdaTn}) we have
$$\nabla\lambda_n=\Big(-\frac{x_2}{ax_1^2},\frac{1}{ax_1}\Big)$$
and then
\begin{equation}\label{nablaun}
|\nabla\lambda_n|=\sqrt{\frac{x_2^2}{a^2x_1^4}+\frac{1}{a^2x_1^2}}=\sqrt{\frac{x_2^2+x_1^2}{a^2x_1^4}}=\frac{1}{|x_1|}\sqrt{\frac{x_2^2}{a^2x_1^2}+\frac{1}{a^2}}\leq\frac{1}{|x_1|}\sqrt{1+\frac{1}{a^2}} 
\end{equation}
Moreover, since $u(A)=0$, applying Morrey's inequality, we obtain that:
\begin{equation}\label{mor}
|u(x)|=|u(x)-u(A)|\leq C||\nabla u||_{r,T_n}(x_1^2+x_2^2)^{\frac{1}{2}-\frac{1}{r}}\leq C||\nabla u||_{r,T_n}(x_1^2+a^2x_1^2)^{\frac{1}{2}-\frac{1}{r}}= C||\nabla u||_{r,T_n}(1+a^2)^{\frac{1}{2}-\frac{1}{r}}|x_1|^{1-\frac{2}{r}}.
\end{equation}
So, by (\ref{nablaun}), (\ref{mor}) and non-negativity of $x_1$, we get
$$\int_{T^n}|(\nabla\lambda_n)u|^sdx\leq C_a||\nabla u||^s_{r,T_n}\int_0^{\frac{1}{2\cdot 3^n}}\Big(\int_0^{ax_1}x_1^{-\frac{2s}{r}}\dx_2\Big)\dx_1=aC_a||\nabla u||^s_{r,T_n}\int_0^{\frac{1}{2\cdot 3^n}}x_1^{1-\frac{2s}{r}}\dx_1=$$
\begin{equation}\label{inte}
=\frac{aC_a}{2-\frac{2s}{r}}||\nabla u||^s_{r,T_n}\Big(\frac{1}{2\cdot 3^n}\Big)^{2-\frac{2s}{r}}=\overline{C}_a||\nabla u||_{r,T_n}^s\cdot (3^n)^{-2+\frac{2s}{r}},
\end{equation}
with 
$$\overline{C}_a=\frac{C}{(2-\frac{2s}{r})\cdot2^{2-\frac{2s}{r}}}a\Big(1+\frac{1}{a^2}\Big)^{\frac{s}{2}}(1+a^2)^{\frac{s}{2}-\frac{s}{r}}.$$
Thus, putting together (\ref{norm}) and (\ref{inte}), we obtain
$$||w_n||^s_{1,s}\leq2^{s-1}||u||^s_{1,s}+C(3^n)^{-2+\frac{2s}{r}}\sum_{i=1}^{3\cdot4^n}\Big(\int_{T^n_i}|\nabla u|^r\dx\Big)^{\frac{s}{r}},$$
with $C=2^{s-1}\overline{C}_a$.\\
Now, applying  H\"{o}lder inequality for sums, with conjugate exponents $\frac{r}{s}$ and $\frac{r}{r-s}$, we have
$$||w_n||^s_{1,s}\leq2^{s-1}||u||^s_{1,s}+C(3^n)^{-2+\frac{2s}{r}}\Big(\sum_{i=1}^{3\cdot4^n}\int_{T^n_i}|\nabla u|^r\dx\Big)^{\frac{s}{r}}\cdot (3\cdot4^n)^{1-\frac{s}{r}}=$$
$$=2^{s-1}||u||^s_{1,s}+3^{1-\frac{s}{r}}C(3^n)^{-2+\frac{2s}{r}}\Big(\sum_{i=1}^{3\cdot4^n}\int_{T^n_i}|\nabla u|^r\dx\Big)^{\frac{s}{r}}\cdot (3^n)^{d_f(1-\frac{s}{r})}$$
\begin{equation}\label{disfin}
=2^{s-1}||u||^s_{1,s}+3^{1-\frac{s}{r}}C\Big(\int_{\Sigma_-^n}|\nabla u|^r\dx\Big)^{\frac{s}{r}}\cdot (3^n)^{-2+\frac{2s}{r}+d_f-\frac{s}{r}d_f}.
\end{equation}
By observing the second term in the last member of the previous chain, we have that it goes to $0$, as $n\to\infty$; in fact:
$$\int_{\Sigma_-^n}|\nabla u|^rdx\to0, \text{ as } n\to\infty, \text{ since }|\Sigma_-^n|\to0$$
and 
$$-2+\frac{2s}{r}+d_f-\frac{s}{r}d_f=-2(1-\frac{s}{r})+d_f(1-\frac{s}{r})=-(2-d_f)(1-\frac{s}{r})<0.$$
\noindent Thus $w_n\in W_0^{1,s}(\Omega_\alpha)$, $\forall 1<s<r$ and $||u||_{1,s,\Omega_\alpha}\leq C$, with $C$ independent on $n$.\\

\noindent To show \textit{(iii)}, we have that  
$$||w_n-u||^s_{1,s,\Omega_\alpha}=||w_n-u||^s_{1,s}\leq C||\nabla(w_n-u)||^s_{s}=C\int_{\Omega_\alpha}|\nabla(w_n-u)|^s\dx=$$
$$=C\int_{\Sigma_-^n}|\nabla(\lambda_nu-u)|^s\dx+C\int_{\Omega_\alpha\setminus\Omega_\alpha^n}|-\nabla u|^s\dx\to0,\text{ as }n\to\infty,$$
since $|\Omega_\alpha\setminus\Omega_\alpha^n|\to0$, as $n\to\infty.$
\endproof

\begin{rem}\label{OpeW2} 
It is possible to prove that:\\
(i) if $u\in C_0(\overline{\Omega}_\alpha)$, then $w_n\in C_0(\overline{\Omega}_\alpha);$\\
(ii) if $u\in W^{1,\infty}_0(\Omega_\alpha)$, then $w_n\in W^{1,\infty}_0(\Omega_\alpha).$\\
Moreover, we observe that $w_n$ in $(ii)$ has a (possibly) different Lipschitz constant with respect to the one of $u$ and it is again independent on $n$.
\end{rem}

\begin{cor}\label{OpeWg}
Let $u$ be in $W_g^{1,r}(\Omega_\alpha)$, $r>2$, with $g\in W^{1,\infty}(\Omega_\alpha)$. Then, the function $z_n(x)=\lambda_n(x)\cdot u(x)+(1-\lambda_n(x))g(x)$, where $\lambda_n(x)$ is defined in (\ref{CoeLam}), has the following properties:\\
\begin{eqnarray*}
	(i) & \,\,\, & z_n(x)\in W_g^{1,s}(\Omega_\alpha), \forall 2<s<r;\\
	(ii) & \,\,\, & ||z_n||_{1,s,\Omega_\alpha}\leq C, \text{ with } C \text{ independent on } n;\\
	(iii) & \,\,\, & z_n\to u \text{ in } W^{1,s}(\Omega_\alpha), \text{ as } n\to\infty.
\end{eqnarray*}
\end{cor}
\proof By its definition, we can write $z_n(x)=g(x)+\lambda_n(x)v(x)$, with $v(x)=u(x)-g(x)$, for each $x\in\Omega_\alpha$.\\
Since $v(x)\in W^{1,p}_0(\Omega_\alpha)$, we get our thesis applying Theorem \ref{OpeW}.
\endproof

\begin{rem}\label{OpeWg2} 
It is possible to prove that:\\
(i) if $u\in C_g(\overline{\Omega}_\alpha)$, then $z_n\in C_g(\overline{\Omega}_\alpha);$\\
(ii) if $u\in W^{1,\infty}_g(\Omega_\alpha)$, then $z_n\in W^{1,\infty}_g(\Omega_\alpha).$\\
Moreover, we observe that $z_n$ in $(ii)$ has a (possibly) different Lipschitz constant with respect to the sum of the ones of $u$ and $g$ and it is again independent on $n$.
\end{rem}

By using the previous result, we now prove Theorem \ref{AsymRes}.

\proof 
Applying the same procedure of Theorem 3.1 in \cite{CF} (for instance) we get that the sequence $\{\tilde{u}_{p,q,n}\}_{n\in\N}$ is bounded in $W_g^{1,p}(\Omega_\alpha)$. Then, there exists $v\in W_g^{1,p}(\Omega_\alpha)$ and a subsequence of $\tilde{u}_{p,q,n}$, that we denote again with $\tilde{u}_{p,q,n}$, such that $\tilde{u}_{p,q,n}\to v$ weakly in $W_g^{1,p}(\Omega_\alpha)$. So, we have
$$J_{p,q}(v)=\frac{1}{p}\int_{\Omega_\alpha}(k^2+|\nabla v|^2)^\frac{p}{2}\dx+\frac{1}{q}\int_{\Omega_\alpha}(k^2+|\nabla v|^2)^\frac{q}{2}\dx-\int_{\Omega_\alpha}fv\dx\leq$$
$$\leq\liminf_{n\to\infty}\Big(\frac{1}{p}\int_{\Omega_\alpha}(k^2+|\nabla\tilde{u}_{p,q,n}|^2)^\frac{p}{2}\dx+\frac{1}{q}\int_{\Omega_\alpha}(k^2+|\nabla\tilde{u}_{p,q,n}|^2)^\frac{q}{2}\dx-\int_{\Omega_\alpha}f_n\tilde{u}_{p,q,n}\dx\Big)=$$
$$=\liminf_{n\to\infty}\Big(\frac{1}{p}\int_{\Omega^n_\alpha}(k^2+|\nabla u_{p,q,n}|^2)^\frac{p}{2}\dx+\frac{1}{q}\int_{\Omega^n_\alpha}(k^2+|\nabla u_{p,q,n}|^2)^\frac{q}{2}\dx-\int_{\Omega^n_\alpha}f_nu_{p,q,n}\dx+$$
$$+\frac{1}{p}\int_{\Omega_\alpha\setminus\Omega^n_\alpha}(k^2+|\nabla g|^2)^\frac{p}{2}\dx+\frac{1}{q}\int_{\Omega_\alpha\setminus\Omega^n_\alpha}(k^2+|\nabla g|^2)^\frac{q}{2}\dx-\int_{\Omega_\alpha\setminus\Omega^n_\alpha}f_ng\dx\Big)\leq$$
$$\leq\limsup_{n\to\infty}\Big(\frac{1}{p}\int_{\Omega^n_\alpha}(k^2+|\nabla u_{p,q,n}|^2)^\frac{p}{2}\dx+\frac{1}{q}\int_{\Omega^n_\alpha}(k^2+|\nabla u_{p,q,n}|^2)^\frac{q}{2}\dx-\int_{\Omega^n_\alpha}f_nu_{p,q,n}\dx+$$
$$+\frac{1}{p}\int_{\Omega_\alpha\setminus\Omega^n_\alpha}(k^2+|\nabla g|^2)^\frac{p}{2}\dx+\frac{1}{q}\int_{\Omega_\alpha\setminus\Omega^n_\alpha}(k^2+|\nabla g|^2)^\frac{q}{2}\dx-\int_{\Omega_\alpha\setminus\Omega^n_\alpha}f_ng\dx\Big)\leq$$
$$\leq\limsup_{n\to\infty}J_{p,q,n}(u_{p,q,n})+$$
$$+\limsup_{n\to\infty}\Big(\frac{1}{p}\int_{\Omega_\alpha\setminus\Omega^n_\alpha}(k^2+|\nabla g|^2)^\frac{p}{2}\dx+\frac{1}{q}\int_{\Omega_\alpha\setminus\Omega^n_\alpha}(k^2+|\nabla g|^2)^\frac{q}{2}\dx-\int_{\Omega_\alpha\setminus\Omega^n_\alpha}f_ng\dx\Big)=$$
\begin{equation}\label{b1}
=\limsup_{n\to\infty}\min_{w\in\mathcal{H}_{p,n}}J_{p,q,n}(w)
\end{equation}
Since, $u_{p,q}$ is the unique solution to Problem (\ref{Ppq}), if we show that $J_{p,q}(v)=J_{p,q}(u_{p,q})$, we will get our thesis.\\
Now, thanks to the Theorem \ref{SumRes}, for $x\in\Omega_\alpha$, let us consider the functions
$$v_n(x)=(w_n(x)\vee\varphi_{1,n}(x))\land\varphi_{2,n}(x), \text{ with } w_n(x)=\lambda_n(x)u_{p,q}(x)+(1-\lambda_n(x))g(x)$$
 and let us show that:\\
$$ \text{(a) } v_n\in\mathcal{H}_{p,n};$$
$$ \text{(b) } v_n\to u_{p,q} \text{ strongly in } W^{1,p}(\Omega_\alpha).$$

Let us prove (a).\\
$v_n\in W_g^{1,p}(\Omega^n_\alpha)$ by Corollary \ref{OpeWg} and the fact that $\varphi_{1,n}\leq g\leq\varphi_{2,n}$ on $\partial\Omega_\alpha^n$. Finally, the fact that $\varphi_{1,n}\leq v_n\leq\varphi_{2,n}$ follows by the definition of $v_n$.\\

 Now, let us prove (b).\\
It follows by \text{(iii)} of Corollary \ref{OpeWg}, the fact that $\varphi_{i,n}\to\varphi_i\in W^{1,p}(\Omega_\alpha)$, for $i=1,2$, and the fact that $\varphi_1\leq u_{p,q}\leq\varphi_2$ in $\Omega_\alpha$. 

Hence, we have
\begin{equation}\label{b2}
\limsup_{n\to\infty}\min_{w\in\mathcal{H}_{p,n}}J_{p,n}(w)\leq \limsup_{n\to\infty}J_{p,q,n}(v_n)=J_{p,q}(u_{p,q}).
\end{equation}
Thus, by $(\ref{b1})$, $(\ref{b2})$ and the fact that $J_p(u_{p,q})\leq J_p(v)$, we get that $v=u_{p,q}$ and then the whole sequence $\tilde{u}_{p,q,n}$ converge to $u_{p,q}$. Furthermore, we obtain that
$$J(u_{p,q})=\lim_{n\to\infty}J_{p,q,n}(u_{p,q,n})$$
and the proof is over.
\endproof

 We can also perform the asympotic analysis for   $u_{\infty,q,n}$ solutions to Problems (\ref{Pqn}) or (\ref{PqLn})    when $n$ goes to $\infty.$ 
 
 In this case the convergence result is achieved by considering the  problems separately:
more precisely for the  solutions to Problems (\ref{Pqn}) we can apply the convergence result of  Theorem 4.1 of  \cite{CF}
and  for the  solutions to Problems (\ref{PqLn}) we can apply the convergence result of  Theorem 4.2 of  \cite{CF}.

Now, let us consider $u_{\infty,q,n}$ solutions to Problems (\ref{Pqn}) or (\ref{PqLn})  respectively and define
\begin{equation}\label{Est2}
\tilde{u}_{\infty,q,n}(x):=\begin{cases}
u_{\infty,q,n}(x), & x\in\overline{\Omega}^n_{\alpha}\\
g(x), & x\in\overline{\Omega}_{\alpha}\setminus\Omega^n_{\alpha}\,\, .
\end{cases}
\end{equation}

\begin{thm}

Let $f_n,f\in L^{p'}(\Omega_\alpha)$, $g\in W^{1,\infty}(\Omega_\alpha)$, $\varphi_{i,n},\varphi_i\in W^{1,p}(\Omega_\alpha)$, for $i=1,2$. Moreover let us assume $\mathcal{H}_{p,n}\neq\emptyset$, $\mathcal{H}_{p}\neq\emptyset$ and, as $n\to\infty$,
$$f_n\to f \text{ in } L^1(\Omega_\alpha)\,\,\, \text{ and }\,\,\, \varphi_{i,n}\to\varphi_{i},\,\, i=1,2,\,\, \text{ in } W^{1,p}(\Omega_\alpha).$$ 
Then $\tilde{u}_{\infty,q,n}(x)$ defined in (\ref{Est2}) admit a subsequence which $*$-weakly converge to solution to Problem  (\ref{Pq1}) or  (\ref{PqL1}) respectively.
\end{thm}

\begin{rem} We briefly discuss about uniqueness. Beside being interesting in itself,  it is a crucial issue in order to obtain the possibility to switch the order of the limits with respect to $n$ and $p$. \\
In \cite{F}, uniqueness results for p-Laplacian unilateral problems are stated (see also \cite{ACJ}, \cite{BDM}, \cite{J} and the references quoted there for the problem of the uniqueness).\\
In our situation, we point out that for the case $L^2+k^2>1$ the issue of the uniqueness is still an open problem both for the fractal and pre-fractal case.\end{rem}

\begin{rem}
We point out that it is possible to extend the present result to other domains with prefractal  and fractal boundaries  like, for example, quasi-filling fractal layers or random snowflakes
(see \cite{CDO} and the reference therein); the key tool is that the domains have  good   \lq\lq extension" properties (see \cite{J}). 
Moreover, it  is possible to perform asymptotic analysis also in the so-called \lq\lq Sobolev admissible domains" (see \cite{Dek}, \cite{Hinz}).
We remark that it is also possible to consider  these  problems  on  fractals structures like, for example, the Sierpinski gasket, where a notion of infinity harmonic functions has been introduced recently (see \cite{CCV}).

\end{rem}

\subsection*{Acknowledgment}

The corresponding author is member of GNAMPA(INdAM)
and is partially supported by Grants Ateneo \lq\lq Sapienza" 2022.

\end{document}